\documentclass[a4paper,12pt, twoside, pdftex]{amsart}

\usepackage{fullpage}
\usepackage{amsmath, amsthm, amssymb}
\usepackage{txfonts}
\usepackage{amsfonts}
\usepackage{prettyref}
\usepackage{mathrsfs}
\usepackage[all]{xy}
\usepackage[utf8]{inputenc}

\usepackage{tikz}
\usetikzlibrary{backgrounds}

\usepackage[OT2,T1]{fontenc}
\DeclareSymbolFont{cyrletters}{OT2}{wncyr}{m}{n}
\DeclareMathSymbol{\Sha}{\mathalpha}{cyrletters}{"58}

\newtheorem{dummy}{dummy}[section]
\newtheorem{lemma}[dummy]{Lemma}
\newtheorem{theorem}[dummy]{Theorem}

\newtheorem{conjecture}[dummy]{Conjecture}
\newtheorem{corollary}[dummy]{Corollary}
\newtheorem{proposition}[dummy]{Proposition}
\theoremstyle{definition}
\newtheorem{definition}[dummy]{Definition}
\newtheorem*{definition*}{Definition}

\newtheorem{remark}[dummy]{Remark}

\newtheorem{assumption}[dummy]{Assumption}

\newtheorem*{acknowledgements}{Acknowledgements}
\newtheorem*{notations}{Notations and conventions}

\numberwithin{equation}{section}

\newrefformat{th}{Theorem~\ref{#1}}
\newrefformat{cr}{Corollary~\ref{#1}}
\newrefformat{lm}{Lemma~\ref{#1}}
\newrefformat{dl}{Definition-Lemma~\ref{#1}}
\newrefformat{df}{Definition~\ref{#1}}
\newrefformat{cl}{Claim~\ref{#1}}
\newrefformat{sl}{Sublemma~\ref{#1}}
\newrefformat{pr}{Proposition~\ref{#1}}
\newrefformat{cj}{Conjecture~\ref{#1}}
\newrefformat{st}{Step~\ref{#1}}
\newrefformat{sc}{Section~\ref{#1}}
\newrefformat{df}{Definition~\ref{#1}}
\newrefformat{rm}{Remark~\ref{#1}}
\newrefformat{q}{Question~\ref{#1}}
\newrefformat{pb}{Problem~\ref{#1}}
\newrefformat{cd}{Condition~\ref{#1}}
\newrefformat{eg}{Example~\ref{#1}}
\newrefformat{he}{Heore~\ref{#1}}
\newrefformat{fg}{Figure~\ref{#1}}
\newrefformat{tb}{Table~\ref{#1}}
\newrefformat{as}{Assumption~\ref{#1}}

\newcommand{\pref}{\prettyref}

\newcommand{\Aut}{\operatorname{Aut}}
\newcommand{\Bl}{\operatorname{Bl}}
\newcommand{\Br}{\operatorname{Br}}

\newcommand{\Coker}{\operatorname{Coker}}

\newcommand{\ddiv}{\operatorname{div}}
\newcommand{\Def}{\operatorname{Def}}
\newcommand{\End}{\operatorname{End}}

\newcommand{\Gr}{\operatorname{Gr}}
\newcommand{\Hom}{\operatorname{Hom}}
\newcommand{\id}{\operatorname{id}}

\newcommand{\Isom}{\operatorname{Isom}}
\newcommand{\JJoin}{\operatorname{Join}}
\newcommand{\Ker}{\operatorname{Ker}}

\newcommand{\ord}{\operatorname{ord}}
\newcommand{\perf}{\operatorname{perf}}

\newcommand{\Pic}{\operatorname{Pic}}
\newcommand{\pr}{\operatorname{pr}}
\newcommand{\pt}{\operatorname{pt}}

\newcommand{\rank}{\operatorname{rank}}

\newcommand{\Spec}{\operatorname{Spec}}

\newcommand{\supp}{\operatorname{supp}}

\newcommand{\tot}{\operatorname{tot}}

\newcommand{\Twist}{\operatorname{Twist}}

\newcommand{\WC}{\operatorname{WC}}

\newcommand{\cP}{\mathcal{P}}

\newcommand{\cS}{\mathcal{S}}

\newcommand{\cU}{\mathcal{U}}

\newcommand{\cX}{\mathcal{X}}
\newcommand{\cY}{\mathcal{Y}}

\newcommand{\bA}{\mathbb{A}}
\newcommand{\bC}{\mathbb{C}}

\newcommand{\bP}{\mathbb{P}}
\newcommand{\bQ}{\mathbb{Q}}

\newcommand{\bZ}{\mathbb{Z}}

\newcommand{\bff}{\mathbf{f}}
\newcommand{\bfg}{\mathbf{g}}

\newcommand{\bfG}{\mathbf{G}}

\newcommand{\bfP}{\mathbf{P}}

\newcommand{\bfS}{\mathbf{S}}
\newcommand{\bfU}{\mathbf{U}}

\newcommand{\bfX}{\mathbf{X}}
\newcommand{\bfY}{\mathbf{Y}}

\newcommand{\scrA}{\mathscr{A}}

\newcommand{\scrE}{\mathscr{E}}
\newcommand{\scrF}{\mathscr{F}}

\newcommand{\scrK}{\mathscr{K}}
\newcommand{\scrL}{\mathscr{L}}
\newcommand{\scrM}{\mathscr{M}}
\newcommand{\scrN}{\mathscr{N}}
\newcommand{\scrO}{\mathscr{O}}

\newcommand{\scrT}{\mathscr{T}}

\newcommand{\scrU}{\mathscr{U}}
\newcommand{\scrV}{\mathscr{V}}

\newcommand{\frakI}{\mathfrak{I}}

\title{Reduced Tate--Shafarevich group}
\author[H.~Morimura]{Hayato Morimura}
\address{Kavli Institute for the Physics and Mathematics of the Universe (WPI),
University of Tokyo,
5-1-5 Kashiwanoha,
Kashiwa,
Chiba,
277-8583,
Japan.}
\email{hayato.morimura@ipmu.jp}

\date{}
\begin{document}
\maketitle

\begin{abstract}
We prove a sort of reconstruction theorem for generic elliptic Calabi--Yau $3$-folds in the sense of C\u{a}ld\u{a}raru.
From our argument it follows that
two generic elliptic Calabi--Yau $3$-folds are derived-equivalent linear over the base
if and only if
their generic fibers are derived-equivalent.
As an application,
we give affirmative answers to the conjectures raised in
\cite{KSS}.
Namely,
for each pair of elliptic Calabi--Yau $3$-folds in
\cite[Table 19]{KSS}
we prove that
they share the relative Jacobian
and
are $\bP^2$-linear derived-equivalent.
\end{abstract}

\section{Introduction}
Among Calabi--Yau manifolds,
there are two classes of considerable interest for both
algebraic geometers
and
string theorists.
One consists of
\emph{Fourier--Mukai partners},
pairs of nonbirational derived-equivalent Calabi--Yau $3$-folds. 
The other consists of
\emph{elliptic Calabi--Yau manifolds},
those which admit elliptic fibrations.
In this paper,
we provide new examples of derived-equivalent elliptic Calabi--Yau $3$-folds.
Our method can be applied to other suitable settings to find more examples.
In fact,
the construction
\cite{KSS}
due to Knapp--Scheidegger--Schimannek is
a special case of orbital degeneracy loci
\cite{BFMT}
and
closed to that of determinantal varieties
\cite{JKLMR}.
We will demonstrate how to apply our method to Inoue varieties in Section
$7$.

First,
we prove the following reconstruction theorem for generic elliptic Calabi--Yau $3$-folds in the sense of C\u{a}ld\u{a}raru.
 
\begin{theorem}[Corollary \pref{cor:reconstruction-CYfibrations}] \label{thm:INTRO1}
Let
$f \colon X \to S, g \colon Y \to S$
be generic elliptic Calabi--Yau $3$-folds.
If there exists an $S$-linear exact equivalence
$\Phi \colon D^b(X) \xrightarrow{\sim} D^b(Y)$,
then 
$g$
is an almost coprime twisted power of
$f$
in the sense of Definition
\pref{dfn:coprime}.
\end{theorem}

The opposite direction is
\cite[Theorem 6.1]{Cal02}.
Nontrivial part is the identification of the relative Jacobians of
$f, g$.
If one knew that
$f, g$
share the same relative Jacobian,
then the result would obtained
either by
\cite[Theorem 1.5]{AKW}
or
careful reading of
\cite{Cal}.
Different respectively from
\cite[Theorem 1.5]{AKW}
and
\cite{Cal},
here
$\dim_\bC X = 3$
and
the Calabi--Yau condition
are crucial to identify the relative Jacobians of
$f, g$.

Presumably,
this is the best possible reconstruction result for generic elliptic Calabi--Yau $3$-folds.
Indeed,
if the smooth parts
$f^\circ \colon X^\circ \to S^\circ , g^\circ \colon Y^\circ \to S^\circ$
are represented by
$\alpha^\circ, \beta^\circ \in \Br^\prime(J^\circ)$,
for any analytic small resolution
$\rho \colon \bar{J} \to J$
we have the $S$-linear exact equivalences
\begin{align*}
D^b(X) \simeq D^b(\bar{J}, \bar{\alpha}), \
D^b(Y) \simeq D^b(\bar{J}, \bar{\beta})
\end{align*} 
from
\cite[Theorem 5.1]{Cal02},
where
$\bar{\alpha}, \bar{\beta} \in \Br_{\text{an}}(\bar{J})$
denote the Brauer classes canonically determined by
$\alpha^\circ, \beta^\circ$.
Moreover,
each
$\rho_i \colon \bar{J}_i \to J_i$
over
$U_i$
with only one node can be either of the two possible resolutions,
which is the source of the differences between
$X, Y$
and
$\bar{J}$.
Hence the right hand sides should not recover 
$X \setminus X^\circ, Y \setminus Y^\circ$.

As an application of
\pref{thm:INTRO1},
in Section
$4$,
we developed a deformation method to construct families of pairs of (derived-equivalent) elliptic Calabi--Yau $3$-folds
which are mutually almost coprime twisted powers. 
Due to
\cite{Wil94, Wil98}
and
\cite[Theorem 1.1]{Mor23},
the morphisms
$f, g$
in the above statement deform to families of derived-equivalent elliptic Calabi--Yau $3$-folds
$\bff \colon \bfX \to \bfS, \bfg \colon \bfY \to \bfS$
over a smooth affine $\bC$-variety
$\Spec B$.
We check that
the deformed Fourier--Mukai kernel is supported on
$\bfX \times_\bfS \bfY$.

Recently,
Knapp--Scheidegger--Schimannek constructed
$12$
pairs of Calabi--Yau $3$-folds
admitting elliptic fibrations over
$\bP^2$
with $5$-section
\cite{KSS}.
Although not smooth,
they are flat
and
have no multiple fibers.
Moreover,
all reducible fibers are
isolated
and
of type
$I_2$.  
The idea was to consider fiberwise homological projective duality
\cite{Kuz06}
for Grassmannian
$\Gr(2, V_5)$
of $2$-planes in
$V_5 \cong \bC^5$
and
its dual
$\Gr(2, V^\vee_5)$
with respective Pl\"{u}cker embeddings into
$\bP(\wedge^2 V_5)$
and
$\bP(\wedge^2 V^\vee_5)$.

\begin{definition}
Let
$A_i, B_i, i = 1, \ldots, 12$
be one of the
$12$
pairs of elliptic Calabi--Yau $3$-folds over
$\bP^2$
labeled as
$i_a, i_b$
in
\cite[Table 19]{KSS}.
We call
$A_i, B_i$
\emph{type $i$ KSS varieties}.
We denote by
$f_i, g_i$
the elliptic fibrations
$A_i \to \bP^2, B_i \to \bP^2$
induced by the canonical projections.
\end{definition}

Based on F-theoretical observations,
Knapp--Scheidegger--Schimannek raised

\begin{conjecture}[{\cite{KSS}}] \label{conj:common}
The elliptic fibrations
$f_i, g_i$
share the relative Jacobian
$\pi_i \colon J_i \to \bP^2$.
\end{conjecture}

By construction it is natural to ask

\begin{conjecture}[{\cite{KSS}}] \label{conj:deq}
Type
$i$
KSS varieties
$A_i, B_i$
are derived-equivalent.
\end{conjecture}

For
$i = 11$
the statement is trivial,
as 
$A_{11}, B_{11}$
are isomorphic.
For
$i = 1, 2$
the statement should follow from
\cite[Remarks 2.3.3, 2.4.3]{KSS}
and
\cite[Proposition 3.5]{Ino}.
Explicitly,
$Y_2, Y_1$
and
$X_2, X_1$
in
\cite{Ino}
which admit elliptic fibrations over
$\bP^2$
are claimed to be respectively isomorphic to
$A_1, A_2$
and
$B_1, B_2$.

In this paper,
we give affirmative answers to these conjectures by proving

\begin{theorem} [\pref{thm:main}] \label{thm:INTRO3}
The elliptic fibrations
$f_i, g_i$
are mutually an
\emph{almost coprime twisted power}
of the other in the sense of Definition
\pref{dfn:coprime}.
\end{theorem}

In particular,
the smooth parts
$f^\circ_i, g^\circ_i$
of
$f_i, g_i$
are respectively isomorphic to the relative moduli spaces of stable sheaves of
rank
$1$
with some degrees
$k, l$
on the fibers of
$g^\circ_i, f^\circ_i$.
Here,
$k, l \in \bZ$
are respectively coprime to the fiber degree of
$g^\circ_i, f^\circ_i$.
Then one can apply
\cite[Theorem 5.1, 6.1]{Cal02}
to obtain a $\bP^2$-linear Fourier--Mukai transform.
Following the argument,
one will obtain a particular example of specialization of derived equivalence
\cite{Mora, Morb},
where the derived equivalence propagates to the whole families. 

\begin{theorem} \label{thm:INTRO3.5}
Let
$f \colon X \to S, g \colon Y \to S$
be generic elliptic Calabi--Yau $3$-folds.
Then there exists an $S$-linear exact equivalence
$\Phi \colon D^b(X) \xrightarrow{\sim} D^b(Y)$
if and only if
their generic fibers are derived-equivalent. 
\end{theorem}

Combination with
\pref{thm:INTRO1}
and
\cite[Theorem 6.1]{Cal02}
yields

\begin{corollary}
Let
$f \colon X \to S, g \colon Y \to S$
be generic elliptic Calabi--Yau $3$-folds.
Then their generic fibers are derived-equivalent
if and only if
$g$
is an almost coprime twisted power of
$f$
in the sense of Definition
\pref{dfn:coprime}. 
\end{corollary}

Toward
\pref{thm:INTRO3}
one needs to prove that
the generic fibers
$A_{i, \eta}, B_{i, \eta}$
of
$f_i, g_i$
share the Jacobian
$J_{i, \eta}$.
In our setting,
this follows from
\cite[Lemma 2.4]{AKW},
as the fiberwise homological projective duality induces the derived equivalence
\begin{align*}
D^b(A_{i, \eta}) \simeq D^b(B_{i, \eta}).
\end{align*}
Similarly,
one obtains derived equivalences of general fibers of
$f_i, g_i$,
which implies their being isomorphic
\cite[Theorem 7.4]{LT}.
Then one might seek to invoke instead
\cite[Lemma 5.5]{DG}.
Alternatively,
one would try to apply
\cite[Proposition 4.2.2]{Cal},
since
$A_i, B_i$
are
Calabi--Yau.

However,
their proofs seem incomplete.
The former proof ended showing 
$A_{i, \eta}, B_{i, \eta}$
to be twists of the geometric generic fiber
$J_{A_i, \bar{\eta}} \cong A_{i, \bar{\eta}} \cong B_{i, \bar{\eta}} \cong J_{B_i, \bar{\eta}}$.
In general,
not all twists come from torsors
\cite[Proposition 5.3]{Sil}.
Namely,
the $\bar{K}$-isomorphism
$J_{A_i, \bar{\eta}} \to J_{B_i, \bar{\eta}}$
might be any isomorphism fixing the origin.
The latter proof ended showing the relative Jacobians to be isomorphic to minimal Weierstrass models
$W(\scrL, a, b), W(\scrL, a^\prime, b^\prime)$,
which become isomorphic over the base
if and only if
their generic fibers are isomorphic.
See also Remark
\pref{rmk:counter}.

Realizing KSS varieties as two different geometric phases of non-abelian gauged linear sigma models,
they also raised

\begin{conjecture}[{\cite{KSS}}] \label{conj:bir}
For
$i \neq 11$
type
$i$
KSS varieties
$A_i, B_i$
are nonbirational.
\end{conjecture}

Again,
for
$i = 1, 2$
the statement should follow from
\cite[Theorem 3.6]{Ino}.
If
$f_i, g_i$
were smooth,
then by
\pref{thm:INTRO3}
and
\cite[Table 19]{KSS}
the morphism 
$g_i$
would be a nonisomorphic coprime twisted power
of
$f_i$
in the sense of Definition
\pref{dfn:coprime}.
This would imply that
the function fields of the sources are nonisomorphic.

Some computations in
\cite{KSS}
were carried out assuming

\begin{conjecture}[{\cite{KSS}}] \label{conj:TS}
The Tate--Shafarevich group
$\Sha_{\bP^2}(J_{i, \eta})$
is isomorphic to
$\bZ_5$.
\end{conjecture}

If
$f_i, g_i$
were smooth,
then one could adapt the argument as in
\cite[Example 1.18]{DG}
to prove the conjecture,
since
$\Br^\prime(A_i)$
vanishes.
We wonder if our deformation method could be helpful to solve
Conjecture
\pref{conj:bir}
and
\pref{conj:TS}.

Provided
\cite[Remark 2.3.3, 2.4.3]{KSS}, 
\pref{thm:INTRO3}
holds also for
\emph{type}
$1, 2$
\emph{Inoue varieties}
$X_1, Y_1$
and
$X_2, Y_2$,
the Fourier--Mukai partners admitting elliptic fibrations with $5$-section
constructed by Inoue
\cite{Ino}.
This implies that
$Y_1, Y_2$
are certain parts of the relative moduli spaces of semistable sheaves of rank
$1$
with suitable degrees on the fibers of
$X_1, X_2$.
In particular,
we give a partial affirmative answer to the question raised in
\cite[Remark 2.12]{Ino}.
The answer will be completed
if the remaining
\emph{type}
$3$
\emph{Inoue varieties}
$X_3, Y_3$
have only irreducible fibers except isolated of type 
$I_2$,
which presumably can be checked by the same method as in
\cite[Section 5, 6]{KSS}.
Note that
$Y_3$
is claimed to be isomorphic to the first in
\cite[Miscellaneous examples]{KSS}
and
admits only one elliptic fibration over
$\bP^1 \times \bP^1$
\cite[Remark 3.13]{Ino}.

By construction the derived equivalence
\begin{align*}
\Phi_{\cP_j} \colon D^b(X_j) \xrightarrow{\sim} D^b(Y_j), \
\cP_j \in D^b(X_j \times Y_j), \
j = 1, 2, 3
\end{align*}
of
$X_j, Y_j$
follows from homological projective duality for categorical joins developed in
\cite{KP}. 
If
$\Phi_{\cP_j}$
are $S_j$-linear,
then by
\pref{thm:INTRO1}
we would obtain an alternative proof of
\pref{thm:INTRO3}
for Inoue varieties.
Hence it is interesting to see
whether the Fourier--Mukai kernels
$\cP_j$
are supported on the fiber products
$X_j \times_{S_j} Y_j$
so that
$\Phi_{\cP_j}$
become $S_j$-linear.

\begin{remark}
Throughout the paper,
we invoke many results from
\cite{Cal}.
As mentioned above,
the proof of
\cite[Proposition 4.2.2]{Cal}
seems incomplete.
However,
it was used only once in the proof of
\cite[Theorem 4.5.2]{Cal}
to show that
the $k$-th twisted power
$f^k$
of a smooth elliptic fibration
$f \colon X \to S$
has the same relative Jacobian as
$f$.
There the usage of
\cite[Proposition 4.2.2]{Cal}
is not crucial,
as the claim immediately follows from
\cite[Proposition 4.2.3]{Cal}.
Moreover,
the author explicitly constructed the $S$-isomorphism
$J_X \to J_{X^k}$
there in terms of the cut-and-reglue procedure in
\cite[Section 4.5]{Cal}.
\end{remark}

\begin{notations}
Throughout the paper,
all $\bC$-varieties are integral separated scheme of finite type over
$\bC$.
Via Serre's GAGA theorem
we sometimes go from the
algebraic
to
analytic
categories for proper $\bC$-varieties.
For any morphism
$f \colon X \to S$
of smooth $\bC$-varieties,
we consider the canonical $S$-linear structure on
$D^b(X) \simeq \perf(X)$
given by the action
\begin{align*}
\perf(X) \times \perf(S) \to \perf(X), \
(F, G) \mapsto F \otimes_{\scrO_X} f^*G.
\end{align*}
Then for morphisms
$f \colon X \to S, g \colon Y \to S$
of smooth $\bC$-varieties,
an exact functor
$\Phi \colon D^b(X) \to D^b(Y)$
is $S$-linear
if
$\Phi$
respects the $S$-linear structures,
i.e.,
we have functorial isomorphisms
\begin{align*}
\Phi(F \otimes_{\scrO_X} f^*G)
\cong
\Phi(F) \otimes_{\scrO_Y} g^* G, \
F \in \perf(X), G \in \perf(S).
\end{align*}
\end{notations}

\begin{acknowledgements}
The author would like to thank
Antnella Grassi,
Johanna Knapp,
Burt Totaro,
Emanuel Sheidegger
and
Thorsten Schimannek
for
helpful comments
and
correspondence.
He is furthermore grateful to an anonymous referee for careful reading
and
informing him on the example in Remark
\pref{rmk:counter}.
This work was
mainly supported by
SISSA PhD scholarships in Mathematics
and
partially supported by
JSPS KAKENHI Grant Number JP23KJ0341.
\end{acknowledgements}

\section{Ogg--Shafarevich theory over the complex number field}
In this section,
we collect relevant results on classification of elliptic fibrations without multiple fibers.
Main references are
\cite{DG, Cal, Sil}.

\subsection{The Weil--Ch\^{a}telet group}
Let
$E$
be an elliptic curve over a filed
$K$
of characteristic
$0$.
The Weil--Ch\^{a}telet group is a pointed set equipped with a group structure
which classifies torsors for
$E$.

\begin{definition}[{\cite[Section X.2]{Sil}}]
The
\emph{isomorphism group}
$\Isom(E)$
of
$E$
is the group of $\bar{K}$-isomorphism from
$E$
to
$E$.
The
\emph{automorphism group}
$\Aut(E)$
of
$E$
is the $G_{\bar{K} / K}$-invariant subgroup of 
$\Isom(E)$
whose elements preserve the origin
$O \in E$.
We use the same symbol
$E$
to denote the elliptic curve
and
its
\emph{translation group},
the $G_{\bar{K} / K}$-invariant subgroup of 
$\Isom(E)$
of translations
$\tau_p$
for points
$p \in E$.
\end{definition}

\begin{lemma}[{\cite[Proposition X5.1]{Sil}}]
There is a bijection of pointed sets
\begin{align*}
E \times \Aut(E)
\to
\Isom(E), \
(p, \alpha)
\mapsto
\tau_p \circ \alpha,
\end{align*}
identifying
$\Isom(E)$
with the product
$E \times \Aut(E)$
twisted by the natural action of
$\Aut(E)$
on
$E$.
\end{lemma}

\begin{definition}[{\cite[Section X.2]{Sil}}]
A
\emph{twist}
of
$E$
is a smooth $K$-curve
$C$
which is $\bar{K}$-isomorphic to
$E$.
Two twists
$C, C^\prime$
are
\emph{equivalent}
if they are $K$-isomorphic.
We denote by
$\Twist(E/K)$
the set of equivalence classes of twists of
$E$. 
\end{definition}

\begin{lemma}[{\cite[Theorem X2.2]{Sil}}]
There is a canonical bijection of pointed sets
\begin{align*}
H^1_{\text{\'et}}(G_{\bar{K}/K}, \Isom(E))
\to
\Twist(E/K).
\end{align*}
\end{lemma}

\begin{definition}[{\cite[Section X.3]{Sil}}]
A
\emph{torsor}
for
$E$
is a smooth $K$-curve
$C$
equipped with a simply transitive algebraic group $E$-action defined over
$K$.
It is
\emph{trivial}
if
$C(K) \neq \emptyset$.
Two torsors
$C, C^\prime$
are
\emph{equivalent}
if there is an $E$-equivariant $K$-isomorphism
$\theta \colon C \to C^\prime$.
The
\emph{Weil--Ch\^{a}telet group}
$\WC(E/K)$
is the set of equivalence classes of torsors for
$E$. 
\end{definition}

\begin{theorem}[{\cite[Theorem X3.6]{Sil}}]
There is a canonical bijection of pointed sets
\begin{align*}
\WC(E/K)
\to
H^1_{\text{\'et}}(G_{\bar{K}/K}, E).
\end{align*}
In particular,
the image of
$H^1_{\text{\'et}}(G_{\bar{K}/K}, E)$
under the inclusion induced by
$E \subset \Isom(E)$
gives a natural group structure to
$\WC(E/K) \subset \Twist(E/K)$.
\end{theorem}

\begin{theorem}[{\cite[Proposition X5.3]{Sil}}]
The inclusion
$\Aut(E) \subset \Isom(E)$
induces
\begin{align*}
H^1_{\text{\'et}}(G_{\bar{K}/K}, \Aut(E))
\subset
H^1_{\text{\'et}}(G_{\bar{K}/K}, \Isom(E)).
\end{align*}
Let
$\Twist((E, O)/K)$
be the image of
$H^1_{\text{\'et}}(G_{\bar{K}/K}, \Aut(E))$
regarded as a subset of
$\Twist(E/K)$.
If
$C \in \Twist((E, O)/K)$
then
$C(K) \neq \emptyset$.
Conversely,
if
$E^\prime$
is an elliptic curve over
$K$
which is $\bar{K}$-isomorphic to
$E$,
then
$E^\prime$
represents an element of
$\Twist((E, O)/K)$.
\end{theorem}

\begin{remark}
In general,
$C \in \Twist((E, O)/K)$
is not $K$-isomorphic to
$E$.
By
\cite[Proposition 5.4]{Sil}
the group
$\Twist((E, O)/K)$
is canonically isomorphic to
$K^*/(K^*)^n$
where
$n$
becomes equal to
$2, 4$
or
$6$
depending on the $j$-invariant of
$E$.
Hence the elements of
$K^*/(K^*)^n$
correspond to the twists of
$E$
which do not come from torsors.
\end{remark}

\subsection{The Tate--Shafarevich group}
Let
$S$
be a normal integral excellent scheme.
We denote by
$\eta = \Spec k(S)$
its generic point.
Let
$E$
be an elliptic curve over
$k(S)$.
Assume that
the characteristic of
$k(S)$
is
$0$.
In good cases,
the Tate--Shafarevich group
can be described as a certain quotient of Brauer groups
and
classifies smooth elliptic fibrations.

\begin{definition}[{\cite[Section 1]{DG}}]
The
\emph{Tate--Shafarevich group}
$\Sha_S(E)$
associated with
$S$
and
$E$
is the subset
$\bigcap_{s \in S} \Ker(loc_{\bar{s}}) \subset \WC(J_\eta/k(S))$
for natural specialization maps
\begin{align*}
loc_{\bar{s}}
\colon
\WC(E/k(S))
\to
\WC(E(\bar{s})/k(\bar{s})), \
C
\mapsto
C(\bar{s})
=
C \times_{k(s)} k(\bar{s}). 
\end{align*}
\end{definition}

There is a standard cohomological interpretation of
$\Sha_S(E)$.
Consider the exact sequence
\begin{align*}
0
\to
H^1_{\text{\'et}}(S, \iota_* E)
\to
H^1_{\text{\'et}}(\eta, E)
\to
H^0_{\text{\'et}}(S, R^1 \iota_* E)
\to
H^2_{\text{\'et}}(S, \iota_* E)
\to
H^2_{\text{\'et}}(\eta, E)
\end{align*}
where
$\iota \colon \Spec k(S) \to S$
denotes the canonical morphism.
For any
$s \in S$
we have
$(R^1 \iota_* E)_{\bar{s}}
\cong
H^1_{\text{\'et}}(\eta_{\bar{s}}, E(\bar{s}))$
and
the natural homomorphism
\begin{align*}
H^0_{\text{\'et}}(S, R^1 \iota_* E)
\to
\prod_{s \in S} (R^1 \iota_* E)_{\bar{s}} 
\end{align*}
is injective.
Since the composition
$H^1_{\text{\'et}}(\eta, E)
\to
H^0_{\text{\'et}}(S, R^1 \iota_* E)
\to
(R^1 \iota_* E)_{\bar{s}}$
coincides with
$loc_{\bar{s}}$,
one obtains an exact sequence
\begin{align*}
0
\to
H^1_{\text{\'et}}(S, \iota_* E)
\to
\WC(E/k(S))
\to
\prod_{s \in S} \WC (E(\bar{s})/k(\bar{s}))
\end{align*}
and
$\Sha_S(E) = H^1_{\text{\'et}}(S, \iota_* E)$.

\begin{remark}
An element
$C \in \WC(E/k(S))$
maps to
$0$
in
$\WC(E(\bar{s})/k(\bar{s}))$
if and only if
there exists an irreducible \'etale neighborhood
$U \to S$
of
$s$
such that
$C \times_{k(S)} k(U)$
has a rational point over
$k(U)$.
Indeed,
if the image of
$C$
in
$(R^1 \iota_* E)_{\bar{s}}$
is
$0$,
then there exists an irreducible \'etale neighborhood
$U \to S$
of
$s$
such that
the image of
$C$
in
$H^1_{\text{\'et}}(k(U), E \times_{k(S)} k(U))$
is
$0$.
Hence
$\Sha_S(E)$
consists of equivalence classes of \'etale locally trivial torsors.
\end{remark}

In good cases,
there is another cohomological interpretation of
$\Sha_S(E)$.
Let
$X$
be a scheme.

\begin{definition}[{\cite[Definition 1.1.7]{Cal}}]
The
\emph{cohomological Brauer group}
$\Br^\prime(X)$
of
$X$
is defined as
$H^2_{\text{\'et}}(X, \scrO^*_X)_{tors}$.
The
\emph{Brauer group}
$\Br(X)$
of
$X$
is the group of isomorphism classes of Azumaya algebras on
$X$
modulo equivalence relation.
Here,
the group structure is given by tensor products.
Two Azumaya algebras
$\scrA, \scrA^\prime$
on
$X$
are
\emph{equivalent}
if there exists a locally free sheaves
$\scrE, \scrE^\prime$
satisfying
\begin{align*}
\scrA \otimes \underline{\End}(\scrE)
\cong
\scrA^\prime \otimes \underline{\End}(\scrE^\prime).
\end{align*}
\end{definition}

\begin{theorem}[{\cite[Theorem 1.1.4]{Cal}}]
If
$X$
is a smooth over $\bC$,
then
$H^2_{\text{\'et}}(X, \scrO^*_X)$
is torsion
and
for the associated analytic space
$X^h$
we have
\begin{align*}
\Br^\prime(X)
=
\Br^\prime_{\text{an}}(X^h)
=
H^2_{\text{an}}(X^h, \scrO^*_X)_{tors}
=
H^2_{\text{an}}(X^h, \scrO^*_X)
\end{align*} 
where
$\Br^\prime_{\text{an}}(X^h)$
denotes the cohomological Brauer group for a complex analytic space
$X^h$.
\end{theorem}

\begin{theorem}[{\cite[Tag 0A2J]{SP}, \cite{Jon}}] \label{lem:Br=Br'}
If
$X$
is
quasicompact
or
connected,
then
$\Br(X)$
is torsion
and
there is a canonical injection
$\Br(X) \to \Br^\prime(X)$.
If
$X$
is quasiprojective over
$\bC$,
then
$\Br(X)$
surjects onto
$\Br^\prime(X)$.
In particular,
if
$X$
is smooth quasiprojective over
$\bC$,
then we have
\begin{align*}
\Br(X) = \Br^\prime(X) = H^2_{\text{\'et}}(X, \scrO^*_X).
\end{align*}
\end{theorem}

\begin{definition}[{\cite[Definition 2.1]{DG}, \cite[Definition 6.1.5]{Cal}}]
An
\emph{elliptic fibration}
$f \colon X \to S$
is a projective morphism of $\bC$-schemes
whose generic fiber
$X_\eta$
is a genus
$1$
regular $k(S)$-curve
and
all fibers are geometrically connected.
The
\emph{discriminant locus}
$\Delta_f$
of
$f$
is the closed subset of points
$s \in S$
over which the fibers
$X_s$
are not regular.
We denote by
$f^\circ \colon X^\circ \to S^\circ$
its smooth part,
i.e.,
the restriction of
$f$
over
$S \setminus \Delta_f$.
A fiber
$X_s$
over a point
$s \in S$
is
\emph{multiple}
if
$f$
is not smooth at any
$x \in X_s$.
A
\emph{section}
of
$f$
is a morphism
$\sigma_X \colon S \to X$
satisfying
$f \circ \sigma_X = \id$.
An
\emph{$n$-section}
of
$f$
is a closed subscheme to
which the restriction of
$f$
becomes a finite morphism of degree
$n > 1$.
\end{definition}

\begin{remark}
Each component of a multiple fiber must be either
of dimension more than one
or
of dimension one
and
nonreduced at all points.
If an elliptic fibration
$f \colon X \to S$
has no multiple fibers,
then \'etale locally it admits a section.
\end{remark}

\begin{theorem}[{\cite[Corollary 1.17]{DG}}] \label{thm:TS-Br}
Let
$f \colon X \to S$
be a flat elliptic fibration of smooth $\bC$-varieties with section
whose generic fiber is isomorphic to
$E$.
Let
$S^{(1)}$
be the set of points in
$S$
with
$\dim \scrO_{S, s} = 1$.
If the fiber
$X_s$
is geometrically integral over each
$s \in S^{(1)}$,
then we have
\begin{align*}
\Sha_S(E) \cong \Coker(\Br^\prime(S) \to \Br^\prime(X)) 
\end{align*}
where the map
$\Br^\prime(S) \to \Br^\prime(X)$
is given by the pullback.
\end{theorem}

\subsection{The relative Jacobian}
In classifying smooth elliptic fibrations,
the relative Jacobian plays the same role as the Jacobian classifies torsors.

\begin{definition}[{\cite[Definition 4.1.1, 4.2.1, 4.5.1]{Cal}}]
A
\emph{smooth elliptic fibration}
$f \colon X \to S$  
is a smooth projective morphism of smooth $\bC$-varieties
whose fiber
$X_s$
over any point
$s \in S$
is a genus
$1$
regular $k(s)$-curve. 
Fix an $f$-ample line bundle
$\scrO_{X/S}(1)$.
The
\emph{relative Jacobian}  
$\pi \colon J \to S$
is the relative moduli space of stable sheaves of rank
$1$,
degree
$0$
on the fibers of
$f$.
The
\emph{$k$-th twisted power}
$f^k \colon X^k \to S$
of
$f$  
is the relative moduli space of stable sheaves of rank
$1$,
degree
$k \in \bZ$
on the fibers of
$f$.
\end{definition}

\begin{remark} \label{rmk:twisted}
There is a bijective correspondence between
smooth elliptic fibrations
$f \colon X \to S$
with relative Jacobian
$\pi \colon J \to S$
and
elements of
$\Sha_S(J_\eta)$.
Note that
$\pi$
admits a section.
Applying
\pref{lem:Br=Br'}
and
\pref{thm:TS-Br},
one obtains
\begin{align*}
\Sha_{S}(J_\eta)
\cong
\Br^\prime(J) /\Br^\prime(S)
=
\Br(J) /\Br(S). 
\end{align*}
By
\cite[Proposition 4.2.3]{Cal}
\'etale locally
$f$
admits a section
and
there exists an \'etale cover
$\{ U_i \}$
of
$S$
such that
the base changes
$f_i \colon X_i \to U_i$
are isomorphic to
$\pi_i \colon J_i \to U_i$.
Moreover,
$f$
is obtained by gluing
$\pi_i$
along fiberwise translations
$\tau_{ij}
\colon
J_j |_{J_i \times_S J_j}
\to
J_i |_{J_i \times_S J_j}$,
which correspond to degree
$0$
line bundles
$\scrO(\sigma_{J_j} - \sigma_{J_i})$
on
$J_i \times_S J_j$.
Via the interpretation of cohomological Brauer classes as gerbes
\cite[Section 1.1]{Cal},
the collection
$\{ \scrO(\sigma_{J_j} - \sigma_{J_i}) \}$
defines an element
$\alpha$
of
$\Br^\prime(J) = \Br(J)$.
As explained in
\cite[Section 4.3]{Cal},
it is the obstruction against the existence of a universal sheaf on
$X \times_S J$.
In particular,
the local universal sheaves
$\scrU_i$
on
$X_i \times_{U_i} J_i$,
together with isomorphisms on double intersections,
form a $\pr^*_2 \alpha$-twisted sheaf on
$X \times_S J$.
Moreover,
$\alpha$
maps to the element of
$\Sha_S(J_\eta)$
which corresponds to
$f$
\cite[Theorem 4.4.1]{Cal}.
\end{remark}

Taking twisted powers is compatible with the above correspondence in the following sense.

\begin{theorem}[{\cite[Theorem 4.5.2]{Cal}}] \label{thm:4.5.2}
Let
$f \colon X \to S$
be a smooth elliptic fibration.
Fix an $f$-ample line bundle
$\scrO_{X/S}(1)$.
Then any $k$-th twisted power
$f^k \colon X^k \to S$
of
$f$
is a smooth elliptic fibration
which has the same relative Jacobian
$\pi \colon J \to S$
as
$f$.
If
$\alpha \in \Sha_S(J_\eta)$
is an element representing
$f$,
then
$f^k$
is represented by
$\alpha^k$.
\end{theorem}

More generally,
the relative Jacobian
and
twisted powers
are defined as follows.

\begin{definition}[{\cite[Definition 6.4.1]{Cal}}]
Let
$f \colon X \to S$
be a flat elliptic fibration of $\bC$-varieties.
Fix an $f$-ample line bundle
$\scrO_{X/S}(1)$
and
a closed point
$s \in S$.
Let
$P$
be the Hilbert polynomial of
$\scrO_{X_s}$
on
$X_s$
with respect to the polarization given by
$\scrO_{X/S}(1) |_{X_s}$.
Consider the relative moduli space
$M_{X/S}(P) \to S$
of semistable sheaves of Hilbert polynomial
$P$
on the fibers of
$f$.
By the universal property of
$M_{X/S}(P) \to S$
there exists a natural section
$S \to M_{X/S}(P)$
which sends
$s$
to the point
$[\scrO_{X_s}]$
representing
$\scrO_{X_s}$.
Let
$J$
be the unique component of
$M$
which contains the image of this section.
The
\emph{relative Jacobian}
is the restriction
$\pi \colon J \to S$
of the morphism
$M_{X/S}(P) \to S$
to
$J$.
\end{definition}

\begin{definition}[{\cite[Notation 6.6.3]{Cal}}]
Let
$f \colon X \to S$
be a flat elliptic fibration with $n$-section without multiple fibers.
Assume that
all reducible fibers of
$f$
are
isolated
and
of type
$I_2$.
Fix an $f$-ample line bundle
$\scrO_{X/S}(1)$.
Let
$M^k_{X/S}(P) \to S$
be the relative moduli space of semistable sheaves of
rank
$1$,
degree
$k$
on the fibers of
$f$.
Let
$X^k$  
be the union of the components of
$M^k_{X/S}(P)$
which contains a point corresponding to a stable line bundle on a fiber of
$f$.
The
\emph{$k$-th twisted power}
of
$f$
is the restriction
$f^k \colon X^k \to S$
of the morphism
$M^k_{X/S}(P) \to S$
to
$X^k$.
\end{definition}

\begin{remark}
The relative Jacobian
$\pi \colon J \to S$
is a flat elliptic fibration with section
whose discriminant locus
$\Delta_\pi$
equals
$\Delta_f$.
If
$S$
is smooth,
then the restriction over the complement
$S^\circ = S \setminus \Delta_f$
coincides with the relative Jacobian
$\pi^\circ \colon J^\circ \to S^\circ$
for a smooth elliptic fibration
$f^\circ$.
Similarly,
the restriction over
$S^\circ$
of
$f^k$
coincides with the $k$-th twisted power
$f^{\circ k} \colon X^{\circ k} \to S^\circ$.
\end{remark}

\subsection{Minimal Weierstrass fibrations}
Given any elliptic fibration
$f \colon X \to S$
of smooth $\bC$-schemes admitting a section,
one obtains a Weierstrass fibration by contracting all but distinguished components of reducible fibers.
Among birational Weierstrass fibrations over the same base,
the one with the smallest discriminant locus is called minimal.

\begin{definition}
Let
$S$
be a $\bC$-scheme
and
$\scrL$
a line bundle on
$S$.
Take global sections
$a \in H^0(S, \scrL^{\otimes 4}),
b \in H^0(S, \scrL^{\otimes 6})$
such that
$4a^3 + 27 b^2 \in H^0(S, \scrL^{\otimes 12})$
is nonzero.
Consider the projective bundle
\begin{align*}
\pi_{\scrL, a, b}
\colon
\bfP_\scrL
=
\bP_S(\scrO_S \oplus \scrL^{\otimes -2} \oplus \scrL^{\otimes -3})
\to
S.
\end{align*}
We denote by
$\scrO_{\bfP_\scrL / S}(1)$
the line bundle corresponding to the relative hyperplane class.
Let
$W(\scrL, a, b) \subset \bfP_\scrL$
be the closed subscheme defined by the equation
$v^2 w = u^3 + a u w^2 + b w^3$,
where
$u, v$
and
$w$
are respectively given by the global sections of
$\scrO_{\bfP_\scrL / S}(1) \otimes \pi^*_\scrL \scrL^{\otimes 2}$,
$\scrO_{\bfP_\scrL / S}(1) \otimes \pi^*_\scrL \scrL^{\otimes 3}$
and
$\scrO_{\bfP_\scrL / S}(1)$
corresponding to the natural injections of
$\scrL^{\otimes -2}$,
$\scrL^{\otimes -3}$
and
$\scrO_S$
into
$\scrO_S \oplus \scrL^{\otimes -2} \oplus \scrL^{\otimes -3}$.
The canonical projection induces a flat elliptic fibration
$\pi_{\scrL, a, b} \colon W(\scrL, a, b) \to S$
admitting a section
$\sigma_{\scrL, a, b} \colon S \to W(\scrL, a, b)$.
We call
$\pi_{\scrL, a, b}$
the 
\emph{Weierstrass fibration}
associated with
$\scrL, a, b$. 
\end{definition}

\begin{remark} \label{rmk:functorial}
All fibers of
$\pi_{\scrL, a, b}$
are irreducible plane cubic curves.
The discriminant locus of
$\pi_{\scrL, a, b}$
is the support of the Cartier divisor defined by
$4a^3 + 27 b^2$.
The construction of a Weierstrass fibrations is functorial.
Namely,
we have
\begin{align*}
W(h^* \scrL, h^*(a), h^*(b))
\cong
W(\scrL, a, b) \times_S S^\prime
\end{align*}
for any morphism
$h \colon S^\prime \to S$
of $\bC$-schemes.
If
$S$
is smooth,
then
$\sigma_{\scrL, a, b}(S)$
lies in the smooth locus of
$W(\scrL, a, b)$.
\end{remark}

\begin{lemma}[{\cite[Proposition 2.4]{DG}}] \label{lem:model}
Let
$E$
be an elliptic curve over
$K$
with a rational point
$\xi \in E(K)$.
For any smooth $\bC$-scheme
$S$
with
$k(S) \cong K$
there exists a Weierstrass fibration
$\pi_{\scrL, a, b} \colon W(\scrL, a, b) \to S$
whose generic fiber is isomorphic to
$E$.
The closure of
$\xi$
is
$\sigma_{\scrL, a, b}(S)$.
\end{lemma}

\begin{lemma}[{\cite[Theorem 2.1]{Nak}, \cite[Theorem 2.3]{DG}}] \label{lem:Weierstrass}
Let
$f \colon X \to S$
be an elliptic fibration of smooth $\bC$-varieties admitting a section
$\sigma_X \colon S \to X$.
Then there exists a birational $S$-morphism from
$X$
to
$W(\scrL, a, b)$,
which is the contraction of all components of fibers
not intersecting
$\sigma_X(S)$.
Moreover,
$\scrL$
is isomorphic to all of
$\sigma^*_X(\Omega_{X/S}),
f_* \omega_{X/S},
(R^1 f_* \scrO_X)^{-1}$
and
$\scrO_{\sigma_X(S)}(-\sigma_X(S))$
when they are invertible.
\end{lemma}

\begin{definition}[{\cite[Definition 2.6]{DG}}] \label{dfn:minimal}
A Weierstrass fibration
$\pi_{\scrL, a, b}  \colon W(\scrL, a, b) \to S$
is
\emph{minimal}
if there is no effective divisor
$D$
such that
$\ddiv(a) \geq 4D, \ddiv(b) \geq 6D$.
\end{definition}

\begin{remark}
Every Weierstrass fibration is birational to a minimal Weierstrass fibration
\cite[Proposition 2.5]{DG}.
The discriminant locus contains that of the minimal Weierstrass fibration.
\end{remark}

\begin{definition}[{\cite[Definition 2.11]{DG}}] \label{dfn:2.11}
A projective morphism
$f \colon X \to S$
is
\emph{relatively minimal}  
if
$X$
is $\bQ$-factorial
and
has only terminal singularities,
and
if
$K_X.C \geq 0$
for
the canonical divisor
$K_X$
and
any irreducible curve
$C \subset X$
mapping to a point.
\end{definition}

\begin{lemma}[{\cite[Proposition 2.16]{DG}}] \label{lem:minimal}
Let
$f \colon X \to S$
be a relatively minimal elliptic fibration admitting a section.
Assume that
$f_* \omega_{X/S}$
is invertible.
Then the Weierstrass fibration
$\pi_{f_* \omega_{X/S}, a, b} \colon W(f_* \omega_{X/S}, a, b) \to S$
is minimal.
\end{lemma}

\begin{proposition}[{\cite[Proposition 2.17]{DG}}]
Let
$f \colon X \to S$
be a relatively minimal elliptic fibration
and
$\pi_{\scrL, a, b} \colon W(\scrL, a, b) \to S$
a minimal Weierstrass fibration
whose generic fiber is isomorphic to
$E$.
If the Jacobian of
$X_\eta$
is isomorphic to 
$E$,
then
$\Delta_f, \Delta_{\pi_{\scrL, a, b}}$
coincide.
\end{proposition}

\section{Reconstruction of Generic elliptic Calabi--Yau $3$-folds}
In this section,
we prove
\pref{thm:INTRO1}.

\subsection{Generic elliptic $3$-folds}
Categorically,
the cut-and-reglue procedure also works for flat elliptic $3$-folds without multiple fibers
whose only reducible fibers are
isolated
and
of type
$I_2$.

\begin{definition}[{\cite[Definition 6.1.6]{Cal}}] \label{dfn:generic} 
A
\emph{generic elliptic $3$-fold}
$f \colon X \to S$
is an elliptic fibration of smooth $\bC$-varieties
$X, S$
with
$\dim_\bC X = 3$
satisfying:
\begin{itemize}
\item[(1)]
$f$
is flat.
\item[(2)]
$f$
does not have multiple fibers.
\item[(3)]
$f$
admits a multisection.
(In particular,
$f$
does not admit a section.)

\item[(4)]
The discriminant locus
$(\Delta_f)_{red}$
with reduced induced scheme structure is an integral curve in
$S$
with only nodes and cusps as singularities.
\item[(5)]
The fiber over a general point of
$\Delta_f$
is a rational curve with one node.
\end{itemize}
We call
$f$
a
\emph{generic elliptic Calabi--Yau $3$-fold}
if in addition
$X$
is a Calabi--Yau in the strict sence,
i.e.,
we have
$\omega_X \cong \scrO_X$
and
$H^1(X, \scrO_X) = 0$.
\end{definition}

\begin{lemma}[{\cite[Theorem 6.1.9]{Cal}}] \label{lem:classification}
Let
$f \colon X \to S$
be a generic elliptic $3$-fold.
Then over any closed point
$s \in S$
the fiber
$X_s$
is one of the following:
\begin{itemize}
\item
a smooth elliptic curve
when
$s \in S \setminus \Delta_f$;
\item
a rational curve with one node
when
$s \in \Delta_f$
is a smooth point;
\item
a reducible curve of type
$I_2$,
i.e.,
two copies of
$\bP^1$
intersecting transversely at two points
when
$s \in \Delta_f$
is a node;
\item
a rational curve with one cusp
when
$s \in \Delta_f$
is a cusp.
\end{itemize}
Moreover,
for each component
$C$
of all type
$I_2$
fibers,
we have
\begin{align*}
\scrN_{C/X}
\cong
\scrO_C(-1) \oplus \scrO_C(-1).
\end{align*}
\end{lemma}

Applying
\pref{lem:Br=Br'}
and
\pref{thm:TS-Br}
to the smooth part
$\pi^\circ$
of the relative Jacobian
$\pi \colon J \to S$,
one obtains
\begin{align*}
\Sha_{S^\circ}(J_\eta)
\cong
\Br^\prime(J^\circ) /\Br^\prime(S^\circ)
=
\Br(J^\circ) /\Br(S^\circ).
\end{align*}
By the standard purity theorem for the cohomological Brauer group,
we also have
\begin{align*}
\Br^\prime(J^\circ) /\Br^\prime(S^\circ)
=
\Br^\prime(J) /\Br^\prime(S). 
\end{align*}

\begin{definition} \label{dfn:coprime}
Let
$f \colon X \to S$
be a generic elliptic fibration with relative Jacobian
$\pi \colon J \to S$
and
$\alpha^\circ \in \Br^\prime(J^\circ)$
a representative for
$f^\circ \in \Sha_{S^\circ}(J_\eta)$.
We call a generic elliptic fibration
$g \colon Y \to S$
a
\emph{coprime twisted power}
of
$f$
if
$g$
is isomorphic to
$f^k$
for some
$k \in \bZ$
coprime to the order
$\ord([\alpha^\circ])$
in
$\Sha_{S^\circ}(J_\eta)$.
We call
$g \colon Y \to S$
an
\emph{almost coprime twisted power}
of
$f$
if
$g^\circ$
is isomorphic to
$f^{\circ k}$
for some
$k \in \bZ$
coprime to
$\ord([\alpha^\circ])$
and
there exists an analytic open cover
$\{ U_i \}$
of
$S$
such that
$Y \times_S U_i, X^k \times_S U_i$
are $U_i$-isomorphic as an analytic space.
\end{definition}

Passing to the analytic category,
$J$
always admits a small resolution of its singularities.
Our motivation for almost coprime twisted powers comes from the following facts.

\begin{lemma}[{\cite[Theorem 6.4.6]{Cal}}] \label{lem:6.4.6}
For any analytic small resolution of singularities
\begin{align*}
\rho \colon \bar{J} \to J
\end{align*}
of
$J$,
there exists an analytic open cover
$\{ U_i \}$
of
$S$
such that
$X_i = X \times_S U_i, \bar{J}_i = \bar{J} \times_S U_i$
are $U_i$-isomorphic as an analytic space.
\end{lemma}
\begin{proof}
Since
$f$
has no multiple fibers,
analytic locally it admits a section
\cite[Theorem 6.1.8]{Cal}.
Moreover,
as all its reducible fibers are
isolated
and
of type
$I_2$,
over sufficiently small analytic open subset
$U_i$
there is at most one type
$I_2$
fiber of
$f$.
For each component
$C \cong \bP^1$
of a type
$I_2$
fiber,
the normal bundle
$\scrN_{C / X}$
is isomorphic to
$\scrO_C(-1) \oplus \scrO_C(-1)$
by
\pref{lem:classification}.
It is well known that
contraction of such a curve as
$C$
in a $3$-fold yields an ordinary double point.
Hence
$\rho$
resolves the ordinary double points.

By
\cite[Proposition 6.4.2, Theorem 6.4.3]{Cal}
there exist sheaves
$\scrU_i$
on
$X_i \times_{U_i} X_i$
flat over the second factors,
which by universality of
$\pi_i \colon J_i \to U_i$
induce surjective $U_i$-morphisms
$X_i \to J_i$.
These are at most contraction of one of the two components of a type
$I_2$
fiber
which does not intersect the local section.
Note that
such morphisms coincide with the morphisms from
\pref{lem:Weierstrass}.
Hence
$X_i \to J_i$
might differ from
$\rho_i \colon \bar{J}_i \to J_i$
up to whether the contracted component intersect the local sections.
Switching components amounts to performing flops to 
$X_i \to J_i$.
\end{proof}

\begin{lemma}[{\cite[Theorem 6.5.1]{Cal}}] \label{lem:6.5.1}
There exists a unique extension
\begin{align*}
\alpha
\in
\Br(J)
\subset
\Br^\prime(J)
\end{align*}
of
$\alpha^\circ \in \Br^\prime(J^\circ)$.
Let
$\bar{\alpha} = \rho^* \alpha$
for any analytic small resolution
$\rho \colon \bar{J} \to J$
of singularities of
$J$.
Then there exists a $\pr^*_2 \bar{\alpha}$-twisted sheaf
$\bar{\scrU}$
on
$X \times_S \bar{J}$
restricting to the $\pr^*_2 \alpha^\circ$-twisted sheaf on
$X^\circ \times_{S^\circ} J^\circ$,
where
$\alpha^\circ$
is the obstruction for the local universal sheaves against gluing to yield a universal sheaf.
\end{lemma} 
\begin{proof}
Let
$\varphi_i \colon X_i \to \bar{J}_i$
be the $U_i$-isomorphisms of analytic spaces from
\pref{lem:6.4.6}.
We write
$\bar{\scrU}_i$
for the pullbacks of
$\scrU_i$
by
$\id \times_{U_i} \varphi^{-1}_i$.
By construction
$\bar{\scrU}_i$
become local universal sheaves away from type $I_2$ fibers.
Moreover,
the obstruction for
$\bar{\scrU}_i$
against gluing to yield a universal sheaf is unique in a suitable sense
\cite[Theorem 3.3.2]{Cal}. 
Since any double intersection
\begin{align*}
\bar{J}_i \cap \bar{J}_j
=
J_i \cap J_j
\end{align*}
for
$i \neq j$
do not contain type
$I_2$
fibers,
the obstruction defines a unique extension
$\alpha \in H^2_{an}(J, \scrO^*_J)$.
Its restriction to
$J^\circ$
must coincides with
$\alpha^\circ$
by uniqueness of the obstruction.
The fact
$\alpha \in \Br(J)$
follows from
\cite[Remark 6.5.2]{Cal}.
\end{proof}

\begin{theorem}[{\cite[Theorem 5.1]{Cal02}}] \label{thm:5.1}
The relative integral transforms
\begin{align} \label{eq:5.1}
\bar{\Phi}_{\bar{\scrU}}
=
\Phi_{\iota_{S *} \bar{\scrU}}
\colon
D^b(X) \xrightarrow{\sim} D^b(\bar{J}, \bar{\alpha})
\end{align}
with kernel
$\bar{\scrU}$
is an equivalence,
where
$\iota_S \colon X \times_S \bar{J} \hookrightarrow X \times \bar{J}$
denotes the closed immersion.
\end{theorem}

\subsection{Reconstruction of fibrations}
Provided the common relative Jacobian,
it is easy to reconstruct generic elliptic $3$-folds in a suitable sense.

\begin{lemma} \label{lem:appliedCGF}
Let
$f \colon X \to S, g \colon Y \to S$
be flat locally projective morphisms of $\bC$-varieties.
Then any $S$-linear exact equivalence
$\Phi \colon D^b(X) \xrightarrow{\sim} D^b(Y)$
induces a $k(S)$-linear exact equivalence
$\Phi_{k(S)} \colon D^b(X_\eta) \xrightarrow{\sim} D^b(Y_\eta)$
on the generic fibers
$X_\eta, Y_\eta$.
\end{lemma}
\begin{proof}
By
\cite[Proposition 2.15]{HLS}
we may assume that
$f, g$
are smooth over
$S = \Spec R$.
Then one can apply
\cite[Theorem 1.1]{Mora}
to obtain $k(S)$-linear exact equivalences
\begin{align*}
D^b(X_\eta) \xrightarrow{\sim} D^b(X) / D^b_0(X), \
D^b(Y_\eta) \xrightarrow{\sim} D^b(Y) / D^b_0(Y).
\end{align*}
to the Verdier quotients by the full $S$-linear triangulated subcategories spanned by complexes with coherent $R$-torsion cohomology.
Since
$\Phi$
is $S$-linear,
universality of Verdier quotients induces the desired equivalence.
\end{proof}

\begin{corollary} \label{cor:isomorphism-generic}
Let
$f \colon X \to S, g \colon Y \to S$
be generic elliptic $3$-folds.
If there exists an $S$-linear exact equivalence
$\Phi \colon D^b(X) \xrightarrow{\sim} D^b(Y)$,
then their generic fibers
$X_\eta, Y_\eta$
share the Jacobian
$J_\eta$.
\end{corollary}
\begin{proof}
The generic fibers
$X_\eta, Y_\eta$
are derived-equivalent by
\pref{lem:appliedCGF}.
Now,
the claim follows  from
\cite[Lemma 2.4]{AKW}.
\end{proof}

\begin{proposition}
Let
$f \colon X \to S, g \colon Y \to S$
be generic elliptic $3$-folds.
Assume that
$S$
is proper over
$\bC$.
If there exists an $S$-linear exact equivalence
$\Phi \colon D^b(X) \xrightarrow{\sim} D^b(Y)$,
then over any closed point
$s \in S$
the fibers
$X_s, Y_s$
are isomorphic.
In particular,
the discriminant loci
$\Delta_f, \Delta_g$
coincide.
\end{proposition}
\begin{proof}
Since
$X, Y$
are smooth proper over
$\bC$,
the equivalence
$\Phi$
is of Fourier--Mukai type by
\cite{Ola}.
Then 
its $S$-linearity
and
\cite[Proposition 2.15]{HLS}
implies the derived equivalence of
$X_s, Y_s$.
Now,
the claim follows from
\pref{lem:classification}
and
\cite[Theorem 7.4(2)]{LT}.
\end{proof}

\begin{lemma} \label{lem:rec-smooth}
Let
$f \colon X \to S, g \colon Y \to S$
be smooth elliptic $3$-folds with relative Jacobian
$\pi \colon J \to S$.
If there exists an $S$-linear exact equivalence
$\Phi \colon D^b(X) \xrightarrow{\sim} D^b(Y)$,
then 
$g$
is a coprime twisted power of
$f$.
\end{lemma}
\begin{proof}
Let
$\alpha, \beta \in \Br^\prime(J)$
be representatives of
$\Sha_{S}(J_\eta) \cong \Br^\prime(J) / \Br^\prime(S)$
for
$f, g$.
There is an \'etale cover
$\{ U_i \}$
of
$S$
such that
to each
$U_i$
the base changes
$f_i, g_i$
admit sections
\begin{align*}
\sigma_{X, i} \colon U_i \to X_i = X \times_S U_i, \
\sigma_{Y, i} \colon U_i \to Y_i = Y \times_S U_i
\end{align*}
and
$\alpha, \beta$
give the gluing data for
$\{ J_i \}$
to yield
$X, Y$
respectively.
The injection
$\Br^\prime(J) \to \Br^\prime(J_\eta)$
induced by the base change
$J_\eta \to J$
to
$k(S)$
sends
$\alpha, \beta$
to
$\alpha_\eta, \beta_\eta \in \Br^\prime(J_\eta)$.
Then the \'etale cover
$\{ U_i \times_S k(S) \}$
of
$\Spec k(S)$
represents both
$\alpha_\eta$
and
$\beta_\eta$.

Since by
assumption
and
\pref{lem:appliedCGF}
the generic fibers
$X_\eta, Y_\eta$
are derived-equivalent,
one can apply
\cite[Lemma 2.4, Theorem 2.5]{AKW}
to see that
$Y_\eta$
is isomorphic to the moduli space of degree
$k$
line bundles on
$X_\eta$
for some
$k \in \bZ$
coprime to the order
$\ord([\alpha])$.
Hence both
$\alpha^k_\eta$
and
$\beta_\eta$
yield
$g_\eta$
via gluing the \'etale cover
$\{ J_i \times_{U_i} k(U_i) \}$
of
$J_\eta$.
Since the map
$\Br^\prime(J) \to \Br^\prime(J_\eta)$
is injective,
both
$\alpha^k$
and
$\beta$
yield
$g$
via gluing the \'etale cover
$\{ J_i \}$
of
$J$.
By
\pref{thm:4.5.2}
the Brauer class
$\alpha^k$
represents
$f^k$.
Hence
$g$
is isomorphic to
$f^k$.
\end{proof}

\begin{remark}
One could have used
\cite[Theorem 1.5]{AKW}
to obtain
$\beta = \alpha^k$,
provided the inherited Fourier--Mukai transform
$D^b(J_\eta, \alpha_\eta) \xrightarrow{\sim} D^b(J_\eta, \beta_\eta)$
from
$D^b(J, \alpha) \xrightarrow{\sim} D^b(J, \beta)$
for derived categories of twisted coherent sheaves
\cite{Morb}. 
\end{remark}

\begin{lemma} \label{lem:rec-fibrations}
Let
$f \colon X \to S, g \colon Y \to S$
be generic elliptic $3$-folds with relative Jacobian
$\pi \colon J \to S$.
If there exists an $S$-linear exact equivalence
$\Phi \colon D^b(X) \xrightarrow{\sim} D^b(Y)$,
then 
$g$
is an almost coprime twisted power of
$f$.
\end{lemma}
\begin{proof}
Let
$\alpha^\circ, \beta^\circ \in \Br^\prime(J^\circ)$
be representatives of
$\Sha_{S^\circ}(J_\eta) \cong \Br^\prime(J^\circ) / \Br^\prime(S^\circ)$
for the smooth parts
$f^\circ, g^\circ$.
By
\pref{lem:6.5.1}
there exist unique extensions
$\alpha, \beta \in \Br(J)$
of
$\alpha^\circ, \beta^\circ$.
Let
$\bar{\alpha}, \bar{\beta} \in \Br(\bar{J})$
be their images under the map
$\Br(J) \to \Br(\bar{J})$
given by the pullback along any analytic small resolution
$\rho \colon \bar{J} \to J$
of singularities.
Take an analytic open cover
$\{ U_i \}$
of
$S$
representing both
$\bar{\alpha}$
and
$\bar{\beta}$
such that
each
$U_i$
contains at most one node of
$\Delta_f = \Delta_g$.
Then
$X_i, Y_i$
are isomorphic as an analytic space over
$U_i$
by
\pref{lem:6.4.6}.
By
\cite[Proposition 2.15]{HLS}
we have the induced $S^\circ$-linear exact equivalence
$D^b(X^\circ) \xrightarrow{\sim} D^b(Y^\circ)$.
One can apply
\pref{lem:rec-smooth}
to find some
$k \in \bZ$
coprime to
$\ord([\alpha^\circ])$
such that
$g^\circ$
is isomorphic to
$f^{\circ k}$.
As explained in
\cite[Section 6.6]{Cal},
the base changes
$X^k_i$
is $U_i$-isomorphic to
$X_i$
as an analytic space.
\end{proof}

\begin{remark}
In the above proof,
we have not used the assumption on
$X$
to be $3$-dimensional.
However,
this is crucial for
\pref{lem:classification}
on the classification of the fibers,
which in turn is responsible for the construction of local universal sheaves
$\bar{\scrU}_i$
on
$X_i \times_{U_i} \bar{J}_i$
\cite[Theorem 6.4.2]{Cal}.
Hence without the assumption one might not have the equivalence
\pref{eq:5.1}.
\end{remark}

\subsection{Calabi--Yau case}
The triviality of the canonical bundle allows us to detect common relative Jacobians via derived equivalences linear over the base.

\begin{remark} \label{rmk:S}
If
$X$
is Calabi--Yau,
then
$S$
must be either a
rational
or
Enriques
surface by
\cite[Proposition 2.3]{Gro}.
\end{remark}

\begin{lemma} \label{lem:nontrivial}
Let
$f \colon X \to S$
be a generic elliptic Calabi--Yau $3$-fold.
Then the restriction
$f^\circ_* \omega_{X^\circ/S^\circ}
=
f_* \omega_{X/S} |_{S^\circ}$
of
$f_* \omega_{X/S}$
to
$S^\circ = S \setminus \Delta_f$
is nontrivial.
\end{lemma}
\begin{proof}
Otherwise,
by
\pref{lem:Weierstrass}
the restriction
$\pi^\circ \colon J^\circ \to S^\circ$
must coincide with the Weierstrass fibration associated with
$\scrO_{S^\circ}$.
Then we would obtain
$J^\circ = \bP^2 \times E$
for some elliptic curve over
$\bC$.
By Künneth formula
we have
\begin{align*}
H^2(\bP^2 \times E, \scrO_{\bP^2 \times E}) = 0, \
H^3(\bP^2 \times E, \bZ) \cong \bZ^2.
\end{align*}
Since
$\Br^\prime(J^\circ)$
is torsion,
we would obtain
$\Br^\prime(J^\circ) = 0$.
By the standard purity theorem for the cohomological Brauer group,
we would obtain
$\Br^\prime(J) = 0$.
Let
$J^{st}$
be the stable part of
$J$,
the complement of the singularities,
which by definition is isomorphic to 
$\bar{J}^{st} = \rho^{-1}(J^{st})$.
Let
$X^{st}$
be the complement in
$X$
of the type
$I_2$
fibers.
Then
\cite[Theorem 3.3.4]{Cal}
and
the proof of
\cite[Theorem 6.4.3]{Cal}
imply that
there would exist
a universal sheaf
$\scrU^{st}$
on
$X^{st} \times_S X^{st}$
and
an isomorphism
$X^{st} \to J^{st}$.
Recall that
the section
$\sigma_J \colon S \to J$
of
$\pi$
factorizes through
$J^{st}$.
Hence
$f$
would admit a section,
which contradicts to condition
$(3)$
in Definition
\pref{dfn:generic}.
\end{proof}

\begin{lemma} \label{lem:cover}
Let
$f \colon X \to S$
be a generic elliptic Calabi--Yau $3$-fold.
Then there exists an affine open cover
$\{ U_i \}^3_{i=1}$
of
$S$
such that
the restrictions
$\scrV |_{U_i}$
of any vector bundle
$\scrV$
are trivial.
\end{lemma}
\begin{proof}
Recall from Remark
\pref{rmk:S}
that
$S$
must be either a
rational
or
Enriques
surface.
When
$S$
is rational,
by
\cite{CSSV}
there exists an affine open cover
$\{ U_i \}^3_{i=1}$
of
$S$
such that
each
$U_i$ 
is isomorphic to the complement in
$\bA^2$
of a hyperplane.
It is well known that
any vector bundle on such a variety is trivial. 
When
$S$
is an Enriques surface,
by
\cite{Cos, Ver}
there exists a degree
$4$
finite cover
$S \to \bP^2$.
Define
$\{ U_i \}^3_{i=1}$
as the inverse image of the standard affine cover of
$\bP^2$.
\end{proof}

\begin{lemma} \label{lem:pushforward}
Let
$f \colon X \to S$
be a generic elliptic Calabi--Yau $3$-fold.
Then the relative dualizing sheaf
$\omega_{X/S}$
is an invertible $\scrO_X$-module
and
we have
$f_* \omega_{X/S} \cong \omega^{-1}_S$.
\end{lemma}
\begin{proof}
Since
$f$
is flat
and
its fibers are Gorenstein curves,
by
\cite[Tag 0E6R]{SP}
the relative dualizing sheaf
$\omega_{X/S}$
is an invertible $\scrO_X$-module.
For any point
$s \in S$
we have
\begin{align*}
\dim_{k(s)}H^0(X_s, (\omega_{X/S})_s)
=
\dim_{k(s)}H^0(X_s, \omega_{X_s})
=
\dim_{k(s)}H^0(\Spec k(s), (f |_{X_s})_* \scrO_{X_s})
=
1,
\end{align*}
where the second equality follows from exactness of
$\underline{\Hom}_{\widetilde{k(s)}}(-, \widetilde{k(s)})$
and
\begin{align*}
R (f |_{X_s})_* \omega_{k(s)}
=
R (f |_{X_s})_* (f |_{X_s})^! \widetilde{k(s)}
\cong
R \underline{\Hom}_{\widetilde{k(s)}}(R (f |_{X_s})_* \scrO_{X_s}, \widetilde{k(s)}).
\end{align*}
Hence one can apply
\cite[Corollary III 12.9]{Har}
to see that
$f_* \omega_{X/S}$
is an invertible $\scrO_S$-module.
Consider the restriction of the first fundamental exact sequence
\begin{align} \label{eq:SES1}
0
\to
f^{\circ *} \Omega_{S^\circ}
\to
\Omega_{X^\circ}
\to
\Omega_{X^\circ/S^\circ}
\to
0
\end{align}
to the smooth part
$f^ \circ \colon X^\circ \to S^\circ$.
By
$\omega_X \cong \scrO_X$
and
the adjunction formula
we obtain
$\omega_{X^\circ/S^\circ} \cong f^{\circ *} \omega^{-1}_{S^\circ}$.
From
\cite[Corollary III 12.9]{Har}
and
the projection formula
it follows
$f^\circ_* \omega_{X^\circ/S^\circ} \cong \omega^{-1}_{S^\circ}$.
We will extend this to obtain the desired isomorphism.

Let
$\{ U_i \}$
be the affine open cover of
$S$
constructed in the proof of
\pref{lem:cover}.
Since
$f^\circ_* \omega_{X^\circ/S^\circ}$
is nontrivial by
\pref{lem:nontrivial},
we have
$S^\circ \not\subset U_i$
for all
$1 \leq i \leq 3$.
Equivalently,
the complement
$U^c_i \subset S$
of each
$U_i$,
which is a union of curves,
is not contained in
$\Delta_f$.
From the proof of
\pref{lem:cover}
their intersections
$U^c_i \cap U^c_j$
for any
$1 \leq i \neq j \leq 3$
consist of at most finitely many points.
Hence so does the union 
$T \subset S$
of all
$U^c_i \cap U^c_j$
for any
$i \neq j$.
Suppose that
we could extend the isomorphism
$f^\circ_* \omega_{X^\circ/S^\circ} \cong \omega^{-1}_{S^\circ}$
to
$\tilde{S} = S \setminus T$.
Then it would canonically extend to an isomorphism
$f_* \omega_{X/S} \cong \omega^{-1}_S$,
as
$T \subset S$
has codimension
$2$.
Hence we may assume
\begin{align*}
(U_i \setminus U_j) \cap S^\circ
=
U_i \cap U^c_j \cap S^\circ
\neq \emptyset, \
(U_j \setminus U_i) \cap S^\circ
=
U_j \cap U^c_i \cap S^\circ
\neq \emptyset
\end{align*}
for any
$i \neq j$.

Now,
express the line bundle
$\scrN = f_* \omega_{X/S} \otimes_{\scrO_S} \omega_S$
as a pair of compatible collections
$(\{ \scrO_{U_i} \}_i, \{ \varphi_{ij} \psi_{ij} \})$
where
$\varphi_{ij}, \psi_{ij}$
denote the transition functions of
$f_* \omega_{X/S}, \omega_S$.
Namely,
we have the commutative diagrams
\begin{align*}
\begin{gathered}
\xymatrix{
U_{ij} \times \bC \ar^{\varphi^{-1}_j \psi^{-1}_j}[r] \ar_{\varphi_{ij} \psi_{ij}}[d] & \scrN |_{U_{ij}} \ar@{=}[d] \\
U_{ij} \times \bC & \scrN |_{U_{ij}}. \ar_{\varphi_i \psi_i}[l]
}
\end{gathered}
\end{align*}
In other words,
$\scrN$
can be obtained by gluing
$\{ U_i \times \bC \}$
along
$\{ \varphi_{ij} \psi_{ij} \}$.
On the other hand,
by
$f_* \omega_{X^\circ/S^\circ} \cong \omega^{-1}_{S^\circ}$
and
$\scrN^\circ = \scrN |_{S^\circ} \cong \scrO_{S^\circ}$
we obtain
$\varphi^\circ_{ij}
=
\varphi_{ij} |_{U_{ij} \cap S^\circ}
=
1$.
Note that
by the argument in the previous paragraph
$\scrN^\circ$
encodes all the transition functions
$\varphi_{ij}$.
Indeed,
on
$S^\circ$
one detects the restriction
$\varphi^\circ_{ij}$
of
$\varphi_{ij}$
to some nonempty open subset of each
$U_{ij}$.
Thus from
$\varphi^\circ_{ij} = 1$
we obtain
$\varphi_{ij} = 1$
for all double intersections
$U_{ij}$,
which completes the proof.
\end{proof}

\begin{theorem} \label{thm:coincidence}
Let
$f \colon X \to S, g \colon Y \to S$
be generic elliptic Calabi--Yau $3$-folds.
If there exists an $S$-linear exact equivalence
$\Phi \colon D^b(X) \to D^b(Y)$,
then
$f, g$
share the relative Jacobian
$\pi \colon J \to S$.
\end{theorem}
\begin{proof}
Let
$\pi_X \colon J_X \to S, \pi_Y \colon J_Y \to S$
be the relative Jacobians of
$f, g$.
Take any analytic small resolutions of their singularities
$\rho_X \colon \bar{J}_X \to J_X, \rho_Y \colon \bar{J}_Y \to J_Y$.
Recall that
the only singularities of
$J_X, J_Y$
are isolated ordinary double points.
Their complements
$J^{st}_X, J^{st}_Y$
are called the stable parts of
$J_X, J_Y$
and
we have
\begin{align*}
\bar{J}^{st}_X = \rho^{-1}_X(J^{st}_X) \cong J^{st}_X, \
\bar{J}^{st}_Y = \rho^{-1}_Y(J^{st}_Y) \cong J^{st}_Y.
\end{align*}
Since the sections
$\sigma_{J_X} \colon S \to J_X, \sigma_{J_Y} \colon S \to J_Y$
of
$\pi_X, \pi_Y$
factorizes through
$J^{st}_X, J^{st}_Y$,
they canonically extend to that
\begin{align*}
\bar{\sigma}_X \colon S \to \bar{J}_X, \
\bar{\sigma}_Y \colon S \to \bar{J}_Y
\end{align*}
of 
$\bar{\pi}_X = \pi_X \circ \rho_X, \bar{\pi}_Y = \pi_Y \circ \rho_Y$.

Since by Remark
\pref{rmk:CY}
below we have
$\omega_{\bar{J}_X} \cong \scrO_{\bar{J}_X},
\omega_{\bar{J}_Y} \cong \scrO_{\bar{J}_Y}$,
the compositions
$\bar{\pi}_X, \bar{\pi}_Y$
are relatively minimal in the sense of Definition
\pref{dfn:2.11}.
Applying
\pref{lem:Weierstrass},
we obtain birational morphisms
\begin{align*}
\bar{J}_X \to W(\scrL, a, b), \
\bar{J}_Y \to W(\scrM, c, d)
\end{align*}
to some Weierstrass fibrations,
which are the contractions of all components of fibers not intersecting
$\bar{\sigma}_X(S), \bar{\sigma}_Y(S)$.
From the description of
$\rho_X, \rho_Y$
it follows
\begin{align*}
W(\scrL, a, b) = J_X, \
W(\scrM, c, d) = J_Y.
\end{align*}
As we have
\begin{align*}
\bar{\pi}_{X *} \omega_{\bar{J}_X/S} |_{S^\circ}
\cong
f^\circ_* \omega_{X^\circ/S^\circ}
\cong
\omega_{S^\circ}^{-1}
\cong
g^\circ_* \omega_{Y^\circ/S^\circ}
\cong
\bar{\pi}_{Y *} \omega_{\bar{J}_Y/S} |_{S^\circ},
\end{align*}
the same argument as the proof of
\pref{lem:pushforward}
shows
\begin{align*}
\bar{\pi}_{X *} \omega_{\bar{J}_X/S}
\cong
\omega_{S}^{-1}
\cong
\bar{\pi}_{Y *} \omega_{\bar{J}_Y/S}.
\end{align*}
Then by
\pref{lem:Weierstrass}
we obtain
\begin{align*}
J_X
=
W(\omega_{S}^{-1}, a, b), \
J_Y
=
W(\omega_{S}^{-1}, c, d).
\end{align*}

One can apply
\pref{lem:minimal}
to see that
\begin{align*}
\pi_X \colon W(\omega_{S}^{-1}, a, b) \to S, \
\pi_Y \colon W(\omega_{S}^{-1}, c, d) \to S
\end{align*}
are minimal in the sense of Definition
\pref{dfn:minimal}.
From
\cite[Proposition 2.7]{DG}
we obtain an $S$-isomorphism
\begin{align*}
J_X
=
W(\omega_{S}^{-1}, a, b)
\cong
W(\omega_{S}^{-1}, c, d)
=
J_Y,
\end{align*}
as by Corollary
\pref{cor:isomorphism-generic}
the generic fibers
$J_{X, \eta}, J_{Y, \eta}$
of
$\pi_X, \pi_Y$
are isomorphic.
Indeed,
by Remark
\pref{rmk:functorial}
and
\cite[Proposition 2.7]{DG}
the restrictions
\begin{align*}
\begin{gathered}
W(\omega_{S}^{-1}, a, b) |_{\pi_X(J^{st}_X) \cup \pi_Y(J^{st}_Y)}
\cong
W(\omega_{S}^{-1} |_{\pi_X(J^{st}_X) \cup \pi_Y(J^{st}_Y)} , a |_{\pi_X(J^{st}_X) \cup \pi_Y(J^{st}_Y)}, b |_{\pi_X(J^{st}_X) \cup \pi_Y(J^{st}_Y)}), \\
W(\omega_{S}^{-1}, c, d) |_{\pi_X(J^{st}_X) \cup \pi_Y(J^{st}_Y)}
\cong
W(\omega_{S}^{-1} |_{\pi_X(J^{st}_X) \cup \pi_Y(J^{st}_Y)} , c |_{\pi_X(J^{st}_X) \cup \pi_Y(J^{st}_Y)}, d |_{\pi_X(J^{st}_X) \cup \pi_Y(J^{st}_Y)})
\end{gathered}
\end{align*}
are isomorphic over
$\pi_X(J^{st}_X) \cup \pi_Y(J^{st}_Y)$.
Since
$\pi_X(J^{st}_X) \cup \pi_Y(J^{st}_Y)$
has codimension
$2$,
the restrictions of
$\omega^{-1}_S, a, b, c, d$
to 
$\pi_X(J^{st}_X) \cup \pi_Y(J^{st}_Y)$
canonically extend to
$S$.
\end{proof}

\begin{remark}
The extension technique
which completes the above proof
was used for instance in the proof of
\cite[Proposition 2.4]{DG}.
\end{remark}

\begin{corollary} \label{cor:reconstruction-CYfibrations}
Let
$f \colon X \to S, g \colon Y \to S$
be generic elliptic Calabi--Yau $3$-folds.
If there exists an $S$-linear exact equivalence
$\Phi \colon D^b(X) \to D^b(Y)$,
then 
$g$
is an almost coprime twisted power of
$f$.
\end{corollary}
\begin{proof}
The claim follows immediately from
\pref{lem:rec-fibrations}
and
\pref{thm:coincidence}.
\end{proof}

\begin{remark}
It sufficed to assume either
$X$
or
$Y$
is Calabi--Yau.
Without loss of generality we may assume that
$X$
is Calabi--Yau.
The triviality of
$\omega_Y$
follows from the uniqueness of Serre functors.
The vanishing of
$H^1(Y, \scrO_Y)$
follows from
$H^2(X, \scrO_X) = H^2(Y, \scrO_Y)$
and
Serre duality.
Note that
the $S$-linear exact equivalence is in particular $\bC$-linear exact equivalence
and
hence naturally isomorphic to the Fourier--Mukai transform
$\Phi_P$
with kernel
$P \in D^b(X \times Y)$
unique up to isomorphism
\cite[Section 2]{Orl}.
One obtains the induced isometry
\begin{align*}
H^2(X, \scrO_X)
\oplus
H^4(X, \scrO_X)
\cong
H^2(Y, \scrO_Y)
\oplus
H^4(Y, \scrO_Y)
\end{align*} 
from
\cite[Corollary 3.1.13, 3.1.14]{Cal}.
\end{remark}

\begin{remark} \label{rmk:CY}
If
$X$
is Calabi--Yau,
then we have
$\omega_{\bar{J}_X} \cong \scrO_{\bar{J}_X}$
again by the equivalence
\pref{eq:5.1}
and
the uniqueness of Serre functors.
\end{remark}

\section{Deformations of almost coprime twisted powers}
Let
$f \colon X \to S, g \colon Y \to S$
be generic elliptic Calabi--Yau $3$-folds.
Assume that
there exists an $S$-linear exact equivalence
$\Phi \colon D^b(X) \xrightarrow{\sim} D^b(Y)$.
Then by
\cite[Section 2]{Orl}
and
$S$-linearity
$\Phi$
is naturally isomorphic to a relative Fourier--Mukai transform 
\begin{align*}
\bar{\Phi}_{\bar{P} / S}
=
R p_{S *} (\bar{P} \otimes q^*_S (-))
\colon
D^b(X)
\xrightarrow{\sim}
D^b(Y)
\end{align*}
with kernel
$\bar{P} \in \perf(X \times_S Y)$
where
$q_S \colon X \times_S Y \to X, p_S \colon X \times_S Y \to Y$
denote the projections.
Note that
we have
$\bar{\Phi}_{\bar{P} / S} = \Phi_P$
for the pushforward
$P = \tau_{S *} \bar{P}$
along the closed immersion
$\tau_S \colon X \times_S Y \hookrightarrow X \times Y$. 
By
\cite[Theorem 1.1]{Mor23}
there exists a smooth affine $\bC$-variety
$\Spec B$
over which we have
smooth projective versal deformations
$\bfX, \bfY$
of
$X, Y$
and
a deformation
$\bfP \in \perf(\bfX \times_B \bfY)$
of
$P$
defining a relative Fourier--Mukai transform
\begin{align*}
\bar{\Phi}_{\bfP / B}
=
R p_{B *} (\bfP \otimes q^*_B (-))
\colon
D^b(\bfX)
\xrightarrow{\sim}
D^b(\bfY),
\end{align*}
where
$q_B \colon \bfX \times_B \bfY \to \bfX, p_B \colon \bfX \times_B \bfY \to \bfY$
denote the projections.

\begin{lemma}
Up to taking \'etale neighborhood of
$\Spec B$,
also
the base
$S$
and
morphisms
$f, g$
deform to give families over
$\Spec B$
of elliptic fibrations
$\bff \colon \bfX \to \bfS, \ \bfg \colon \bfY \to \bfS$.
\end{lemma}
\begin{proof}
The claim immediately follows from
\cite{Wil94, Wil98}.
\end{proof}

\begin{lemma} \label{lem:generic} 
Up to shrinking
$\Spec B$,
over any closed point
$b \in \Spec B$
the fibers
$\bff_b \colon \bfX_b \to \bfS_b, \bfg_b \colon \bfY_b \to \bfS_b$
are generic elliptic Calabi--Yau $3$-folds.
\end{lemma}
\begin{proof}
By construction
$\bff_b, \bfg_b$
are elliptic Calabi--Yau $3$-folds.
We check that
general fibers of
$\bff, \bfg$
satisfy conditions
(1), \ldots, (5)
in Definition
\pref{dfn:generic}.
By
construction
and
\cite[Corollary IV11.2.7]{GD66}
after shrinking 
$\Spec B$
the morphism
$\bff, \bfg$
becomes flat.
Since
$\bff_b, \bfg_b$
are flat
and
$\bfX_b, \bfY_b$
are smooth,
$\bff_b, \bfg_b$
do not have multiple fibers.
Condition
(3)
follows from the well known fact that
any elliptic Calabi--Yau manifold admits a multisection.
Conditions
(4), (5)
hold over a sufficiently small open neighborhood of
$0 \in \Spec B$.
\end{proof}

\begin{lemma} \label{lem:key}
The kernel
$\bfP \in \perf(\bfX \times_B \bfY)$
is supported on
$\bfX \times_\bfS \bfY$.
\end{lemma}
\begin{proof}
Let
$R \cong \bC \llbracket t_1, \ldots, t_{\dim_{\bC} \text{H}^1 (X, \scrT_X)} \rrbracket$
be the formal power series ring
which prorepresents the deformation functors
$\Def_X, \Def_Y$.
By
\cite[Theorem III5.4.5]{GD61}
there exist effectivizations
$\cX, \cY$
of universal formal families.
Let
$\{ R_\lambda \}_{\lambda \in \Lambda}$
be the filtered inductive system used to algebrize
$\cX, \cY$.
It is a compatible system of finitely generated $\bC[t_1, \ldots, t_{\dim_{\bC} \text{H}^1 (X, \scrT_X)}]$-subalgebras of
$R$
whose colimit is
$R$.
Let
$\cX_{R_\lambda}, \cY_{R_\lambda}$
be the $R_\lambda$-deformations
of
$X, Y$
used to algebrize
$\cX, \cY$.
Their pullbacks along the canonical homomorphism
$R_\lambda \to R$
are isomorphic to
$\cX, \cY$.
Then
$\bfX, \bfY$
are the pullbacks of
$\cX_{R_\lambda}, \cY_{R_\lambda}$
along some homomorphism
$R_\lambda \to B$.
In summary,
we have the commutative diagram
\begin{align*}
\begin{gathered}
\xymatrix{
X \times_S Y \ar@{^{(}->}[r] \ar_{}[d] & \cX \times_\cS \cY \ar[r]^-{} \ar_{}[d] & \cX_{R_\lambda} \times_{\cS_{R_\lambda}} \cY_{R_\lambda} \ar^{}[d] & \bfX \times_\bfS \bfY \ar[l]^-{} \ar_{}[d] \\
\Spec \bC \ar@{^{(}->}[r] & \Spec R \ar[r]_-{} & \Spec R_\lambda & \Spec B \ar[l]_-{} \\
X \times Y \ar@{^{(}->}[r] \ar_{}[u] & \cX \times_R \cY \ar[r]^-{} \ar_{}[u] & \cX_{R_\lambda} \times_{R_\lambda} \cY_{R_\lambda} \ar^{}[u] & \bfX \times_B \bfY. \ar[l]^-{} \ar_{}[u]
}
\end{gathered}
\end{align*}
Note that
the upper vertical arrows are flat projective,
while the lower vertical arrows are smooth projective for sufficiently large
$\lambda \in \Lambda$.

By
\cite[Proposition 3.6.1]{Lie}
there exists an effectivization
$\cP \in \perf(\cX \times_R \cY)$
of a formal $R$-deformation of
$P$.
Let
$\cP_{R_\lambda} \in \perf(\cX_{R_\lambda} \times_{R_\lambda} \cY_{R_\lambda})$
be the perfect complex used to algebrize
$\cP$
\cite[Proposition 2.2.1]{Lie}.
Its derived pullback
$\cP_{R_\lambda} \otimes^L_{R_\lambda} R$
is isomorphic to
$\cP$.
Then
$\bfP$
is the derived pullback of
$\cP_{R_\lambda}$
along the homomorphism
$R_\lambda \to B$
used to algebrize
$\cX, \cY$.
Regarding
$X \times Y$
as a closed subscheme of
$\cX \times_R \cY$,
by
\cite[Lemma 3.29]{Huy}
we have
\begin{align} \label{eq:supp}
\supp(\cP) \cap (X \times Y)
=
\supp(\cP |_{X \times Y})
=
\supp(P)
=
X \times_S Y
\end{align}
and
$\supp(\cP) \subset \cX \times_R \cY$
is a proper closed subset.
Since the structure morphism
$\cX \times_R \cY \to \Spec R$
is flat proper,
it sends
$\supp(\cP)$
to the closed point
which implies
$\supp(\cP) \subset X \times Y$.
From
\pref{eq:supp}
it follows
$\supp(\cP) \subset X \times_S Y$.
In particular,
the restriction
$\cP |_\cU$
to
$\cU = \cX \times_R \cY \setminus \cX \times_\cS \cY$
is acyclic.
Consider the collection
$\{ \cU_\lambda \}_{\lambda \in \Lambda}$
of complements
$\cU_\lambda = \cX_{R_\lambda} \times_{R_\lambda} \cY_{R_\lambda} \setminus \cX_{R_\lambda} \times_{\cS_{R_\lambda}} \cY_{R_\lambda}$,
which are
flat separated $R_\lambda$-schemes of finite presentation.
For
$\lambda^\prime \in \Lambda$
with
$\lambda^\prime > \lambda$
we have
$\cU_\lambda \cong \cU_{\lambda^\prime} \times_{R_{\lambda^\prime}} R_\lambda$
by construction.
Now,
one can apply
\cite[Proposition 2.2.1]{Lie}
to see that
$\cP_\lambda |_{\cU_\lambda}$
is acyclic for sufficiently large
$\lambda$.
Thus the restriction
$\bfP |_\bfU$
to
$\bfU = \bfX \times_B \bfY \setminus \bfX \times_\bfS \bfY$
becomes acyclic after replacing
$\lambda$
if necessary,
which completes the proof.
\end{proof}

\begin{remark}
While the natural projection
$\supp(P) \to X$
is surjective
\cite[Lemma 6.4]{Huy},
by
$\supp(\cP) \subset X \times Y \subsetneq \cX \times_R \cY$
the natural projection
$\supp(\cP) \to \cX$
cannot be surjective.
This is not a contradiction,
as
\cite[Lemma 6.4]{Huy}
is a statement for $\bC$-varieties.
Indeed,
the proof does not work in our setting as there is no closed point in the complement of
$X$
in
$\cX$.
\end{remark}

\begin{theorem} \label{thm:family}
Up to shrinking
$\Spec B$,
there exists
$k(b) \in \bZ$
for each closed point
$b \in \Spec B$
such that
$\bfg_b$
is an almost coprime $k(b)$-th twisted power of
$\bff_b$.
\end{theorem}
\begin{proof}
The claim immediately follows from
Corollary
\pref{cor:reconstruction-CYfibrations},
\pref{lem:generic}
and
\pref{lem:key}.
\end{proof}

\section{Fourier--Mukai partners from KSS varieties}
In this section,
we prove
\pref{thm:INTRO3}.

\subsection{Grassmannian side}
Let
$\scrF^\vee$
be a globally generated vector bundle of rank
$5$
on
$\bP^2$
and
$\bfG = \Gr_{\bP^2}(2, \scrF)$
the Grassmannian bundle
whose fiber over any point
$x \in \bP^2$
is the Grassmannian
$\Gr(2, \tot(\scrF)_x)$
of $2$-planes in the $k(x)$-vector space
$\tot(\scrF)_x$.
We denote
by
$\scrO_{\bfG / \bP^2}(1)$
the line bundle corresponding to the relative hyperplane class
and 
by
$\pi_\bfG$
the canonical projection.
Let
$\scrE^\vee$
be a globally generated homogeneous vector bundle of rank
$5$
on
$\bfG$
and
$s \in H^0(\bfG, \scrE^\vee)$
a general section.
By the generalized Bertini theorem,
the zero locus
$A = Z(s)$
is a smooth projective $3$-fold.
If in addition
$\omega_{\bfG} \cong \det^{-1}\scrE^\vee$
then
$\omega_A$ 
becomes trivial.
Setting
$\scrF^\vee = F$
and
$\scrE^\vee = \scrO_{\bfG / \bP^2}(1) \otimes \pi^*_\bfG E^\prime$
for
$F, E^\prime$
in 
\cite[Table 2]{KSS},
one obtains Calabi--Yau $3$-folds
$A$.
We will put subscript
$i$
on
$\scrF^\vee, \scrE^\vee, F, E^\prime, \bfG$
and
$A$
to specify which row we are dealing with.

\begin{lemma} \label{lem:strict}
The $3$-fold
$A_i$
is Calabi--Yau in the strict sense,
i.e.,
we have
$H^1(A_i, \scrO_{A_i}) = 0$
in addition to
$\omega_{A_i} \cong \scrO_{A_i}$. 
\end{lemma}
\begin{proof}
Concatenating
Koszul resolution of the ideal sheaf
$\frakI_{A_i}$
of
$A_i$
and
the short exact sequence
$0
\to
\frakI_{A_i}
\to
\scrO_{\bfG_i}
\to
\scrO_{A_i}
\to
0$,
we obtain an exact sequence
\begin{align*}
0
\to
\wedge^5 \scrE_i
\to
\wedge^4 \scrE_i
\to
\cdots
\to
\scrE_i
\to
\scrO_{\bfG_i}
\to
\scrO_{A_i}
\to
0.
\end{align*}
Due to the spectral sequences
\begin{align*}
H^q(\bfG_i, \wedge^p \scrE_i)
\Rightarrow
H^{q-p}(A_i, \scrO_{A_i}),
\end{align*}
it suffices to show the vanishing of
$H^{p+1}(\bfG_i, \wedge^p \scrE_i)$
for
$0 \leq p \leq 5$,
which follows from Leray spectral sequence
\begin{align*}
H^s(\bP^2, \wedge^p E^{\prime \vee}_i \otimes R^r \pi_{\bfG_i *} \scrO_{\bfG_i / \bP^2}(-p))
\Rightarrow
H^{r+s}(\bfG_i, \wedge^p \scrE_i).
\end{align*}
\end{proof}

\subsection{Pfaffian side}
Let
$\scrE^\vee$
be a globally generated vector bundle of rank
$5$
on
$\bP^2$
and
$\bfP = \bP_{\bP^2}(\scrE)$
the projective bundle.
We denote
by
$\scrO_{\bfP / \bP^2}(1)$
the line bundle corresponding to the relative hyperplane class
and
by
$\pi_{\bfP}$
the canonical projection.
Let
$\scrF^\vee$
be a globally generated vector bundle of rank
$5$
on
$\bfP$
and
$\phi \colon \scrF \to \scrF^\vee \otimes \scrO_{\bfP / \bP^2}(1)$
a skew-symmetric morphism corresponding to
$s_\phi \in H^0(\bfP, \wedge^2 \scrF^\vee \otimes \scrO_{\bfP / \bP^2}(1))$.
We denote by
$B$
the first nontrivial degeneracy locus
\begin{align*}
D_3(\phi)
=
\{ x \in \bfP | \rank \phi(x) \leq 3 \}
=
\{ x \in \bfP | \rank \phi(x) \leq 2 \}
=
D_2(\phi)
\end{align*}
of
$\phi$.
Since
$\Pic(\bfP)$
has no torsion,
one can apply the first lemma in
\cite[Section 3]{Oko}
to obtain an exact sequence
\begin{align*}
0
\to
\det \scrF^\vee \otimes \scrO_{\bfP / \bP^2}(-2)
\xrightarrow{(\frac{1}{2} \wedge^2 \phi)^T}
\scrF
\xrightarrow{\phi}
\scrF^\vee \otimes \scrO_{\bfP / \bP^2}(1)
\xrightarrow{\frac{1}{2} \wedge^2 \phi}
\frakI_B \otimes \det \scrF^\vee \otimes \scrO_{\bfP / \bP^2}(3)
\to
0,
\end{align*}
where
$\frakI_B$
denotes the ideal sheaf of
$B$.
From the first proposition in
\cite[Section 3]{Oko},
it follows that
$B$
is a smooth projective $3$-fold,
as
$\wedge^2 \scrF^\vee \otimes \scrO_{\bfP / \bP^2}(1)$
is globally generated.
If in addition
$(\det \scrF)^{\otimes 2} \cong (\det \scrE)(3)$
then
$\omega_B$ 
becomes trivial by the second lemma in
\cite[Section 3]{Oko}.
Setting
$\scrE^\vee = E^\prime$
and
$\scrF^\vee = \pi^*_\bfP F$
for
$E^\prime, F$
in 
\cite[Table 2]{KSS},
one obtains Calabi--Yau $3$-folds
$B$.
We will put subscript
$i$
on
$\scrE^\vee, \scrF^\vee, E^\prime, F, \bfP$
and
$B$
to specify which row we are dealing with.

\begin{lemma}
The $3$-fold
$B_i$
is Calabi--Yau in the strict sense,
i.e.,
we have
$H^1(B_i, \scrO_{B_i}) = 0$
in addition to
$\omega_{B_i} \cong \scrO_{B_i}$. 
\end{lemma}
\begin{proof}
Concatenating the locally free resolution of
$\frakI_{B_i}$
from the second lemma in
\cite[Section 3]{Oko}
and
the short exact sequence
$0
\to
\frakI_{B_i}
\to
\scrO_{\bfP_i}
\to
\scrO_{B_i}
\to
0$,
we obtain
\begin{align*}
0
\to
\scrL_{0, i}
\to
\scrF_{0, i}
\to
\scrF^\vee_{0, i} \otimes \scrL_{0, i}
\to
\scrO_{\bfP_i}
\to
\scrO_{B_i}
\to
0
\end{align*}
where
$\scrL_{0, i}
=
(\det \scrF_i)^{\otimes 2} \otimes \scrO_{\bfP_i / \bP^2}(-5)$
and
$\scrF_{0, i}
=
\scrF_i \otimes \det \scrF_i \otimes \scrO_{\bfP_i / \bP^2}(-3)$.
The vanishing of
the first
and
the second
cohomology of
$\scrL_{0, i},
\scrF_{0, i},
\scrF^\vee_{0, i} \otimes \scrL_{0, i}$
and
$\scrO_{\bfP_i}$
follows from Leray spectral sequence.
\end{proof}

\subsection{Elliptic fibrations over $\bP^2$}
According to
\cite[Section 2.3, 2.4]{KSS},
one can apply the main theorem in
\cite{Ogu}
to see that
$\pi_{\bfG_i}, \pi_{\bfP_i}$
respectively restrict to elliptic fibrations
$f_i \colon A_i \to \bP^2, g_i \colon B_i \to \bP^2$
with $5$-sections.
They
are flat
and
have no multiple fibers.
Moreover,
all reducible fibers of
$f_i, g_i$
are
isolated
and
of type
$I_2$
\cite[Section 5, 6]{KSS}.
According to the well known classification of the fibers of elliptic fibrations,
this implies that
$f_i, g_i$
are generic elliptic Calabi--Yau $3$-folds in the sense of Definition
\pref{dfn:generic}.

\begin{remark}
The existence of type
$I_2$
fibers
\cite[Section 5, 6]{KSS}
implies that
the morphisms
$f_i, g_i$
are not smooth.
In
\cite{KSS}
the authors called
$f_i, g_i$
smooth genus one fibrations,
apparently because
$A_i, B_i$
are smooth.
If this is the case,
then smoothness follows automatically from the above constructions,
despite the comment on usage of Higgs transitions in
\cite[Introduction]{KSS}.
\end{remark}

\begin{lemma} \label{lem:gfibers}
The generic fibers of
$f_i, g_i$
are derived-equivalent.
In particular,
they share the Jacobian
$J_{i, \eta}$.
\end{lemma}
\begin{proof}
The fiber of
$f_i$
over a point
$x \in \bP^2$
is given by
\begin{align*}
\Gr(2, V_5) \times_{\bP(\wedge^2 V_5)} \bP(\tot(E^{\prime \vee}_i)^\perp_x)
\end{align*}
with identifications
\begin{align*}
\tot(F^\vee_i)_x \cong V_5, \
\tot(\scrO_{\bfG_i / \bP^2}(1))_x \otimes \pi^*_{\bfG_i} \tot(E^{\prime \vee}_i)^\perp_x
\cong
\tot(E^{\prime \vee}_i)^\perp_x,
\end{align*}
where
$\tot(E^{\prime \vee}_i)^\perp_x
\subset
\wedge^2 \tot(F^\vee_i)_x$
denotes the orthogonal subspace to a fixed inclusion
$\tot(E^{\prime \vee}_i)_x
\subset
\wedge^2 \tot(F_i)_x$.
Observe from the explicit description of
$s \in H^0(\bfG_i, E_i)$
as in
\cite[Section 2.3]{KSS}
that
$s$
defines a $5$-dimensional quotient $k(x)$-vector space of
$\wedge^2 \tot(F^\vee_i)_x$
whose complement is
$\tot(E^{\prime \vee}_i)^\perp_x$.
Then the fiber of
$g_i$
over general
$x \in \bP^2$
is given by
\begin{align*}
\Gr(2, V^\vee_5) \times_{\bP(\wedge^2 V^\vee_5)} \bP(\tot(E^{\prime \vee}_i)_x).
\end{align*}
Note that
the subvariety of
${\bP(\wedge^2 V^\vee_5)}$
defined by the
$4 \times 4$
Pfaffians of a general
$5 \times 5$
skew-symmetric matrix is isomorphic to
$\Gr(2, V^\vee_5)$.
One can apply
\cite[Theorem 1.1]{Kuz07}
to obtain a derived equivalence of the generic fibers.
See
\cite[Theorem 2.24]{KP}
for the same statement over more general base.
Now,
the claim follows from
\cite[Lemma 2.4]{AKW}.
\end{proof}

\begin{remark}
In
\cite[Section 2.5]{KSS}
the authors claimed that
the above fiberwise orthogonal description globalizes to that of
$A_i$
and
$B_i$.
For their global description,
one needs
$E^{\prime \vee}_i$
to be a subbundle of
$\wedge^2 F_i$
up to twisting by line bundles.
However,
each
$F_i$
is a direct sum of line bundles on
$\bP^2$
and
a subbundle of any line bundle is either
$0$
or
itself.
Then most of
$E^{\prime \vee}_i$
in
\cite[Table 2]{KSS}
cannot be a subbundle of
$\wedge^2 F_i$
no matter how twisted.
\end{remark}

\begin{corollary} \label{cor:fibers}
Over any closed point
$x \in \bP^2$
the fibers of
$f_i, g_i$
are isomorphic.
\end{corollary}
\begin{proof}
By
\pref{thm:main}
below there exists a $\bP^2$-linear Fourier--Mukai transform
\begin{align*}
\Phi_i \colon D^b(A_i) \xrightarrow{\sim} D^b(B_i).
\end{align*}
Then the claim follows from
\cite[Proposition 2.15]{HLS}.
\end{proof} 

\subsection{Common relative Jacobian}
Here is a sufficient condition
which works in our setting. 

\begin{proposition} \label{prop:lift}
Let
$f \colon X \to S, g \colon Y \to S$
be flat elliptic $3$-folds over smooth $\bC$-variety
$S$
without multiple fibers.
Assume that
all reducible fibers of
$f$
are of type
$I_2$.
Assume further that the following conditions hold:
\begin{itemize}
\item[(1)]
The generic fibers of
$f, g$
share the Jacobian
$J_\eta$.
\item[(2)]
There exist possibly analytic small resolutions of singularities
\begin{align*}
\rho_X \colon \bar{J}_X \to J_X, \
\rho_Y \colon \bar{J}_Y \to J_Y
\end{align*}
such that
$\bar{\pi}_X = \pi_X \circ \rho_X, \bar{\pi}_Y = \pi_Y \circ \rho_Y$
give relatively minimal elliptic fibrations
and
$\bar{\pi}_{X *} \omega_{\bar{J}_X/S},
\bar{\pi}_{Y *} \omega_{\bar{J}_Y/S}$
are isomorphic invertible sheaves.
\end{itemize}
Then
$f, g$
share the relative Jacobian
$\pi \colon J \to S$.
\end{proposition}
\begin{proof}
We run the same argument as the proof of
\pref{thm:coincidence}.
The relatively minimal elliptic fibrations
$\bar{\pi}_X, \bar{\pi}_Y$
of smooth $\bC$-varieties admit sections,
as the sections
$\sigma_{J_X}, \sigma_{J_Y}$
of
$\pi_X, \pi_Y$
factorize through the complements
$J^{st}_X, J^{st}_Y$
in
$J_X, J_Y$
of the singular loci.
Since
$\rho_X, \rho_Y$
are small resolutions,
the singular loci consist of finitely many points
whose fibers are unions of curves.
By
\pref{lem:Weierstrass},
\pref{lem:minimal}
and
condition
(2)
there exist birational morphisms
\begin{align*}
\bar{J}_X \to W(\bar{\pi}_{X *} \omega_{\bar{J}_X/S} , a, b), \
\bar{J}_Y \to W(\bar{\pi}_{Y *} \omega_{\bar{J}_Y/S}, c, d)
\end{align*}
to minimal Weierstrass fibrations,
which are the contractions of all components of fibers not intersecting
$\bar{\sigma}_X(S), \bar{\sigma}_Y(S)$.
Note that
by assumption such component over each singular point is a copy of
$\bP^1$.
Hence we obtain
\begin{align*}
J_X
=
W(\bar{\pi}_{X *} \omega_{\bar{J}_X/S} , a, b), \
J_Y
=
W(\bar{\pi}_{Y *} \omega_{\bar{J}_Y/S}, c, d).
\end{align*}
Using the same extension technique as in the proof of
\pref{thm:coincidence},
from conditions
(1), (2)
and
\cite[Proposition 2.7]{DG}
we obtain an $S$-isomorphism
$J_X \cong J_Y$.
\end{proof}

\begin{remark}
When
$\rho_X, \rho_Y$
contract all components of fibers not intersecting
$\bar{\sigma}_X(S), \bar{\sigma}_Y(S)$,
one can remove the assumption on reducible fibers to be of type
$I_2$.
\end{remark}

\begin{remark} \label{rmk:counter}
In earlier version of this paper,
condition
(1)
required only very general fibers to be isomorphic.
However,
the following example informed by an anonymous referee implies that
our original proof was wrong.
This example also implies that
\cite[Lemma 5.5]{DG}
cannot be true.
Consider any elliptic fibration
$f \colon X \to S$.
Suppose that
$S$
admits a nontrivial double covering
$T \to S$.
Let
$Y$
be the quotient of
$X_T = X \times_S T$
by
$\bZ_2$,
where the action is given by
involution on
$T$
and
negation on the fibers.
The generic fiber of
$g \colon Y \to S$,
so called
\emph{quadratic twist},
is not isomorphic to that of
$f$.
On the other hand,
over any closed point
$s \in S$
the fibers of
$f, g$
are isomorphic.
\end{remark}

\begin{corollary} \label{cor:common}
The elliptic fibrations
$f_i, g_i$
share the relative Jacobian
$\pi_i \colon J_i \to \bP^2$.
\end{corollary}
\begin{proof}
It suffices to check that
$f_i, g_i$
satisfy the conditions in Proposition
\pref{prop:lift}.
Condition
(1)
follows from
\pref{lem:gfibers}.
As for condition
(2),
take any analytic small resolutions of singularities
\begin{align*}
\rho_{A_i} \colon \bar{J}_{A_i} \to J_{A_i}, \
\rho_{B_i} \colon \bar{J}_{B_i} \to J_{B_i}.
\end{align*}
By Remark
\pref{rmk:CY}
we have
$\omega_{\bar{J}_{A_i}}
\cong
\scrO_{\bar{J}_{A_i}},
\omega_{\bar{J}_{B_i}}
\cong
\scrO_{\bar{J}_{B_i}}$
and
$\bar{\pi}_{A_i}, \bar{\pi}_{B_i}$
are relatively minimal elliptic fibrations.
Although here the discriminant loci
$\Delta_{f_i}, \Delta_{g_i}$
might be reducible,
one can adapt the proof of
\pref{lem:pushforward}
in a straightforward way to obtain
\begin{align*}
\bar{\pi}_{A_i *} \omega_{\bar{J}_{A_i}/\bP^2}
\cong
\omega^{-1}_{\bP^2}
\cong
\bar{\pi}_{B_i *} \omega_{\bar{J}_{B_i}/\bP^2}.
\end{align*}
\end{proof}

\subsection{Derived equivalence}
\begin{theorem} \label{thm:main}
The elliptic fibrations
$f_i, g_i$
are mutually an
\emph{almost coprime twisted power}
of the other in the sense of Definition
\pref{dfn:coprime}.
\end{theorem}
\begin{proof}
By Corollary
\pref{cor:common}
the restrictions of
$f_i, g_i$
over
$V_i = \bP^2 \setminus \Delta_{\pi_i}$
represents some elements
$\alpha_i, \beta_i \in \Br(J_{i, V_i})$.
We use the same symbol to denote their images under the injection
$\Br^\prime(J_{i, V_i}) \to \Br^\prime(J_{i, \eta})$
given by the pullback along the canonical morphism
$J_{i, \eta} \to J_{i, V_i}$.
Since by
\pref{lem:gfibers}
the generic fibers of
$f_i, g_i$
are derived-equivalent,
one can apply
\cite[Lemma 2.4, Theorem 2.5]{AKW} 
to obtain
$\beta_i = \alpha^k_i$
for some
$k \in \bZ$
coprime to the order
$\ord([\alpha_i])$
in
$\Sha_{V_i}(J_{i, \eta})
\cong
\Br^\prime(J_{i, V_i}) / \Br^\prime(V_i)
=
\Br^\prime(J_{i, V_i})$.
Then
$g_{i, V_i}$
is isomorphic to
$f^k_{i, V_i}$
and
$B_{i, U_j}$
become $U_j$-isomorphic to
$A^k_{i, U_j}$
as an analytic space for some analytic open cover
$\{ U_j \}$
of
$\bP^2$
by
\pref{lem:6.4.6}.
\end{proof}

\begin{corollary} \label{cor:deq}
There exists a $\bP^2$-linear Fourier--Mukai transform
$\Phi_i \colon D^b(A_i) \xrightarrow{\sim} D^b(B_i)$.
\end{corollary}
\begin{proof}
By
\pref{thm:main}
we may assume that
$g_{i, V_i}$
is isomorphic to
$f^k_{i, V_i}$.
Then
$g_{i, V_i}$
represents
$\alpha^k_i \in \Br^\prime(J_{i, V_i})$.
Applying
\pref{thm:5.1},
we obtain
$\bP^2$-linear equivalences
\begin{align*}
D^b(A_i) \xrightarrow{\sim} D^b(\bar{J}_i, \bar{\alpha}_i), \
D^b(B_i) \xrightarrow{\sim} D^b(\bar{J}_i, \bar{\alpha}^k_i).
\end{align*}
Then the claim follows from the $\bP^2$-linear equivalence 
\begin{align*}
D^b(\bar{J}_i, \bar{\alpha}_i)
\xrightarrow{\sim}
D^b(\bar{J}_i, \bar{\alpha}^k_i)
\end{align*}
from
\cite[Theorem 6.1]{Cal02}.
\end{proof}

\begin{remark}
Following the above arguments,
one sees that
the derived equivalence of the generic fibers implies that of generic elliptic Calabi--Yau $3$-folds,
which completes the proof of
\pref{thm:INTRO3.5}.
\end{remark}

\section{Inoue varieties as almost coprime twisted powers}
In this section,
we apply our argument to Inoue varieties to see that
the derived equivalences are linear over the bases.

\subsection{Inoue varieties}
Let
$M_1 = \Gr(2, V_5)$
be a Grassmannian of $2$-planes in
$V_5 \cong \bC^5$
and
$M_2 = \bP_{S_i}(\scrE_i)$
a rank
$r_i$
projective bundle over a del Pezzo surface
$S_i$
satisfying the following conditions:
\begin{itemize}
\item[(i)]
$\scrE^\vee_i$
is globally generated.
\item[(ii)]
$\dim_\bC \varphi_{L_2}(\bP_{S_i}(\scrE_i)) \geq r$
where
$\varphi_{L_2}$
denotes the morphism defined by the line bundle
$\scrO_{M_2 / S_i}(L_2)$
corresponding to the relative hyperplane class
$L_2$
of
$\pi_{\scrE_i, S_i} \colon M_2 \to S_i$.
\item[(iii)]
$\det \scrE_i \cong \omega_{S_i}$.
\end{itemize}
We denote by
$\Sigma_1$
and
$\Sigma_2$
the image of
$M_1$
and
$M_2$
under
the Pl\"{u}cker embedding
and
the morphism defined by the relative hyperplane class
respectively.

\begin{lemma}[{\cite[Proposition 3.1]{Ino}}] \label{lem:CY3}
Let
$\bP_{M_1, M_2}
=
\bP_{M_1 \times M_2}(\scrO(-L_1) \oplus \scrO(-L_2))$
be the resolved join of
$M_1$
and
$M_2$,
where
$L_1$
denotes the Schubert divisor class of
$M_1$
and
$L_2$
is the relative hyperplane class of
$M_2$.
Then a general complete intersection
$X$
of
$r_i + 5$
relative hyperplanes in
$\bP_{M_1, M_2}$
is a Calabi--Yau $3$-fold in the strict sense.
\end{lemma}

\begin{remark}
The image of
$\bP_{M_1, M_2}$
under the morphism
$\varphi_H$
coincides with the projective join
$\JJoin(\Sigma_1, \Sigma_2)$
of
$\Sigma_1$
and
$\Sigma_2$,
where
$H$
denotes the relative hyperplane class of
$\pi_{M_1, M_2} \colon \bP_{M_1, M_2} \to M_1 \times M_2$.
In general,
$\JJoin(\Sigma_1, \Sigma_2)$
is singular along the disjoint union
$\Sigma_1 \sqcup \Sigma_2$.
The morphism
$\varphi_H$
gives a resolution of
$\JJoin(\Sigma_1, \Sigma_2)$.
In particular,
the restriction of
$\varphi_H$
to any enough general complete intersection
$X$
becomes an isomorphism. 
\end{remark}

Let
$\scrE^\perp_i$
be the orthogonal locally free sheaf of
$\scrE_i$.
Namely,
we have a short exact sequence
\begin{align*}
0 \to \scrE^\perp_i \to H^0(S_i, \scrE^\vee_i) \otimes \scrO_{S_i} \to \scrE^\vee_i \to 0.
\end{align*}
We denote by
$r^\prime_i$
and
$L^\prime_2$
the rank of
$\scrE^\perp_i$
and
the relative hyperplane class of
$\pi_{\scrE^\perp_i, S_i} \colon \bP_{S_i}(\scrE^\perp_i) \to S_i$
respectively.
Assume the following additional conditions:
\begin{itemize}
\item[(iv)]
$\dim_{\varphi_{L^\prime_2}}(\bP_{S_i}(\scrE^\perp_i)) \geq r^\prime$.
\item[(v)]
$H^1(S_i, \scrE_i) = 0$.
\end{itemize}
Then
$(\scrE^\perp_i)^\vee$
is globally generated
and
$\det \scrE^\perp_i \cong \omega_{S_i}$.

\begin{corollary} \label{cor:CY3}
Let
$\bP_{M^\prime_1, M^\prime_2}
=
\bP_{M^\prime_1 \times M^\prime_2}(\scrO(-L^\prime_1) \oplus \scrO(-L^\prime_2))$
be the resolved join of
$M^\prime_1 = \Gr(2, V^\vee_5)$
and
$M^\prime_2 = \bP_S(\scrE^\perp_i)$,
where
$L^\prime_1$
denotes the Schubert divisor class of
$M^\prime_1$
and
$L^\prime_2$
is the relative hyperplane class of
$M^\prime_2$.
Then a general complete intersection
$Y$
of
$r^\prime_i + 5$
relative hyperplanes in
$\bP_{M^\prime_1, M^\prime_2}$
is a Calabi--Yau $3$-fold in the strict sense.
\end{corollary}

Consider the cases
where
$M_2 = N_i = \bP_{S_i}(\scrE_i),
M^\prime_2 = N^\prime_i = \bP_{S_i}(\scrE^\perp_i)$
are one of the following:
\begin{itemize}
\item[(1)]
$N_1
=
\bP_{\bP^2} (\scrO_{\bP^2}(-1)^{\oplus^3})
=
\bP^2 \times \bP^2$, \
$N^\prime_1
=
\bP_{\bP^2} (\scrK^{\oplus^3}_1)$,
\item[(2)]
$N_2
=
\bP_{\bP^2} (\scrO_{\bP^2}(-2) \oplus \scrO_{\bP^2}(-1))
=
\Bl_{\pt} \bP^3$, \
$N^\prime_2
=
\bP_{\bP^2} (\scrK_2 \oplus \scrK_1)$,
\item[(3)]
$N_3
=
\bP_{\bP^1 \times \bP^1} (\scrO_{\bP^1 \times \bP^1}(-1, -1)^{\oplus^2})
=
\bP^1 \times \bP^1 \times \bP^1$, \
$N^\prime_3
=
\bP_{\bP^1 \times \bP^1} (\scrK^{\oplus^2}_{1, 1})$.
\end{itemize}
Here,
$\scrK_{1, 1}, \scrK_j$
for
$j = 1, 2$
denote respectively the kernel of the surjections
\begin{align*}
H^0(\bP^1 \times \bP^1, \scrO_{\bP^1 \times \bP^1}(1, 1)) \otimes \scrO_{\bP^1 \times \bP^1}
\to
\scrO_{\bP^1 \times \bP^1}(1, 1), \
H^0(\bP^2, \scrO_{\bP^2}(j)) \otimes \scrO_{\bP^2}
\to
\scrO_{\bP^2}(j).
\end{align*} 
We write
$V_{N_i}$
for
$H^0(N_i, \scrO(L_2))^\vee$.
Let
$W_i \subset \wedge^2 V_5 \oplus V_{N_i}$
be general codimension
$r_i + 5$
linear subspaces
and
$W^\perp_i \subset \wedge^2 V^\vee_5 \oplus V^\vee_{N_i}$
their orthogonal subspaces.
By
\pref{lem:CY3}
and
Corollary
\pref{cor:CY3}
the complete intersections
\begin{align*}
X_i
=
\bP_{M_1, N_i} \times_{\bP(\wedge^2 V_5 \oplus V_{N_i})} \bP(W_i), \
Y_i
=
\bP_{M^\prime_1, N^\prime_i} \times_{\bP(\wedge^2 V^\vee_5 \oplus V^\vee_{N_i})} \bP(W^\perp_i)
\end{align*}
of
$r_i + 5, r^\prime_i + 5$
relative hyperplanes in
$\bP_{M_1, N_i}, \bP_{M^\prime_1, N^\prime_i}$
are Calabi--Yau $3$-folds in the strict sense.
For
$i = 1, 2, 3$
we call
$X_i, Y_i$
\emph{type $i$ Inoue varieties}.
By
\cite[Proposition 3.5, Theorem 3.6]{Ino}
$X_i, Y_i$
are nonbirational derived-equivalent.

\subsection{Elliptic fibrations of Inoue varieties}
Consider the compositions
\begin{align*}
\begin{gathered}
\varpi_{X_1} \colon X_1 \hookrightarrow \bP_{M_1, N_1} \to M_1 \times N_1 \to N_1 = \bP^2 \times \bP^2 \xrightarrow{\pr_1} \bP^2, \\
\varpi_{X_2} \colon X_2 \hookrightarrow \bP_{M_2, N_2} \to M_2 \times N_2 \to N_2 = \bP_{\bP^2}(\scrE_2) \xrightarrow{\pi_{\scrE_2, \bP^2}} \bP^2, \\
\varpi_{X_3} \colon X_3 \hookrightarrow \bP_{M_3, N_3} \to M_3 \times N_3 \to N_3 = \bP^1 \times \bP^1 \times \bP^1 \xrightarrow{\pr_{1, 2}} \bP^1 \times \bP^1  
\end{gathered}
\end{align*}
which are shown to be elliptic fibrations
\cite[Lemma 3.7]{Ino}.
Similarly,
$Y_i$
admit elliptic fibrations
$\varpi_{Y_i}$
over
$\bP^2$
for
$i = 1, 2$
and
$\bP^1 \times \bP^1$
for
$i = 3$
\cite[Lemma 3.11, Remark 3.13]{Ino}.

\begin{lemma} \label{lem:gfibers-Inoue}
The generic fibers of
$\varpi_{X_i}, \varpi_{Y_i}$
for
$i = 1, 2, 3$
are derived-equivalent.
In particular,
they share the Jacobian
$J^\prime_{i, \eta}$.
\end{lemma}
\begin{proof}
The fibers of
$\varpi_{X_1}, \varpi_{Y_1}$
over a point
$x \in \bP^2$
are given by
\begin{align*}
\JJoin(\Gr(2, V_5), \{ x \} \times \bP^2) \cap \bP(W), \
\JJoin(\Gr(2, V^\vee_5), \{ x \} \times \bP^2) \cap \bP(W^\perp).
\end{align*}
Hence general fibers are linear sections of
$\Gr(2, V_5), \Gr(2, V^\vee_5)$
of codimension
$5$.
By definition of
$W, W^\perp$,
they respectively coincide with
\begin{align*}
\Gr(2, V_5) \times_{\bP(\wedge^2 V_5)} \bP(W_x), \
\Gr(2, V^\vee_5) \times_{\bP(\wedge^2 V^\vee_5)} \bP(W^\perp_x)
\end{align*}
for some $5$-dimensional subspace
$W_x \in \wedge^2 V_5$
and
its orthogonal subspace
$W^\perp_x \in \wedge^2 V_5$.
We have similar dual descriptions of the fibers also for
$i = 2, 3$.
Now,
the claim follows from the same argument as in
\pref{lem:gfibers}.
\end{proof}

In the sequel,
we will suppose the following condition for completeness.

\begin{assumption} \label{asp:Inoue}
For
$i = 1, 2$
the elliptic fibrations
$\varpi_{X_i}, \varpi_{Y_i}$
are flat
and
have no multiple fibers.
Moreover,
all reducible fibers are
isolated
and
of type
$I_2$.
\end{assumption}

\begin{remark}
According to
\cite[Remarks 2.3.3, 2.4.3]{KSS}
the Calabi--Yau $3$-folds
$Y_1, Y_2$
and
$X_1, X_2$
are respectively isomorphic to
$A_2, A_1$
and
$B_2, B_1$.
However,
the author could not understand the mathematical reason
why this is the case.
Presumably,
one can directly check the above condition
but we will not pursue it here.
\end{remark}

\begin{corollary} \label{cor:fibers-Inoue}
Under Assumption
\pref{asp:Inoue},
the fibers of
$\varpi_{X_i}, \varpi_{Y_i}$
for
$i = 1, 2$
over any closed point
$x \in \bP^2$
are isomorphic.
\end{corollary}
\begin{proof}
Provided Assumption
\pref{asp:Inoue},
the proof of Corollary
\pref{cor:fibers}
carries over.
Note that
all reducible fibers of
$g_2, g_1$
are
isolated
and
of type
$I_2$.
By
\cite[Lemma 3.11, Remark 3.13]{Ino}
the Calabi--Yau $3$-folds
$X_1, X_2$
admit only one elliptic fibration.
\end{proof}

\subsection{Derived equivalence of Inoue varieties revisited}
\begin{lemma} \label{lem:common-Inoue}
Under Assumption
\pref{asp:Inoue},
type
$i$
Inoue varieties
$X_i, Y_i, i = 1, 2$
share the relative Jacobian
$\varpi_i \colon J^\prime_i \to \bP^2$.
\end{lemma}
\begin{proof}
Provided Assumption
\pref{asp:Inoue},
the proof of Corollary
\pref{cor:common}
carries over,
since
the generic fibers of
$\varpi_{X_i}, \varpi_{Y_i}$
share the Jacobian
and
all reducible fibers of
$\varpi_{X_i}, \varpi_{Y_i}$
are
isolated
and
of type
$I_2$.
\end{proof}

\begin{remark}
According to
\cite[Remark 2.3.3]{KSS}
the Calabi--Yau $3$-fold
$Y_3$
is isomorphic to the first in
\cite[Miscellaneous examples]{KSS}.
Then the elliptic fibration from it must be isomorphic to
$\varpi_{Y_3}$,
as
$Y_3$
admit only one elliptic fibration
\cite[Remark 3.13]{Ino}.
By
\pref{lem:gfibers-Inoue}
the generic fibers of
$\varpi_{X_3}, \varpi_{Y_3}$
share the Jacobian.
If all reducible fibers of
$\varpi_{X_3}, \varpi_{Y_3}$
are
isolated
and
of type
$I_2$,
then from Proposition
\pref{prop:lift}
it follows that
$\varpi_{X_3}, \varpi_{Y_3}$
share the relative Jacobian
$\varpi_3 \colon J^\prime_3 \to \bP^2$.
\end{remark}

\begin{theorem} \label{thm:main-Inoue}
Under Assumption
\pref{asp:Inoue},
the elliptic fibrations
$\varpi_{X_i}, \varpi_{Y_i}$
for
$i = 1, 2$
are mutually an almost coprime twisted power of the other.
\end{theorem}
\begin{proof}
Provided
Assumption
\pref{asp:Inoue}
and
\pref{lem:common-Inoue},
the proof of
\pref{thm:main}
carries over.
\end{proof}

\begin{remark}
Suppose that
all reducible fibers of
$\varpi_{X_3}, \varpi_{Y_3}$
are
isolated
and
of type
$I_2$.
Then one similarly shows
\pref{thm:main-Inoue}
for
$i = 3$,
since we have
$\Br^\prime(\bP^1 \times \bP^1) = 0$
by the standard purity theorem for the cohomological Brauer group,
as
$\bP^1 \times \bP^1$
is rational.
\end{remark}

Now,
from the same arguments as in
Corollary
\pref{cor:deq}
we obtain

\begin{corollary}
Under Assumption
\pref{asp:Inoue},
type
$i$
Inoue varieties
$X_i, Y_i$
for
$i = 1, 2$
are $\bP^2$-linear derived-equivalent.
\end{corollary}


\end{document}